\documentclass[a4paper,12pt]{amsart}
\usepackage{tcolorbox}
\usepackage{pdflscape}
\usepackage[justification=centering,labelformat=empty]{caption}
\usepackage{amsopn,amsfonts,amssymb,amsthm,mathtools,amsmath}
\usepackage{color}
\usepackage{dsfont}
\usepackage{latexsym}

\usepackage{wrapfig}
\usepackage[T1]{fontenc}
\usepackage{lmodern}
\usepackage{textcomp}
\usepackage{soul} 
\usepackage{fancyhdr}

\usepackage{hyperref}
\hypersetup{colorlinks=true,linkcolor=blue,citecolor=blue}
\usepackage[utf8]{inputenc}
\usepackage{float}
\usepackage{hhline}
\usepackage{multirow}
\usepackage{anysize}
\marginsize{2cm}{2cm}{2cm}{2cm}
\usepackage{graphicx}

\def\u{\mathfrak{u}}

\def\g{\mathfrak{g}}
\def\h{\mathfrak{h}}
\def\n{\mathfrak{n}}

\def\s{\mathfrak{s}}
\def\v{\mathfrak{v}}
\def\m{\mathfrak{m}}

\def\e{\operatorname{e}}

\def\u{\mathfrak{u}}
\def\gl{\mathfrak{gl}}

\def\X{\mathfrak{X}}

\def\B{\mathcal{B}}

\renewcommand{\O}{\mathcal{O}}
\def\l{\ell}

\def\C{\mathbb{C}}
\def\R{\mathbb{R}}

\def\Z{\mathbb{Z}}
\def\N{\mathbb{N}}

\def\H{\mathbb{H}}

\newcommand{\GL}{\operatorname{GL}}
\newcommand{\SL}{\operatorname{SL}}
\newcommand{\SU}{\operatorname{SU}}

\def\ad{\operatorname{ad}}

\def\Tr{\operatorname{Tr}}

\def\I{\operatorname{I}}

\newcommand{\Aut}{\operatorname{Aut}}

\newcommand{\Ker}{\operatorname{Ker}}
\newcommand{\Der}{\operatorname{Der}}

\def\alt{\raise1pt\hbox{$\bigwedge$}}
\def\pint{\langle \cdotp,\cdotp \rangle }
\def\la{\langle}
\def\ra{\rangle}
\def\multiset#1#2{\ensuremath{\left(\kern-.3em\left(\genfrac{}{}{0pt}{}{#1}{#2}\right)\kern-.3em\right)}}

\theoremstyle{plain}
\newtheorem{theorem}{\bf Theorem}[section]
\newtheorem{corollary}[theorem]{\bf Corollary}
\newtheorem{proposition}[theorem]{\bf Proposition}
\newtheorem{lemma}[theorem]{\bf Lemma}

\theoremstyle{definition}
\newtheorem{definition}[theorem]{\bf Definition}
\newtheorem{example}[theorem]{\bf Example}

\theoremstyle{remark}
\newtheorem{remark}[theorem]{\bf Remark}

\newcommand{\ri}{{\rm (i)}}
\newcommand{\rii}{{\rm (ii)}}
\newcommand{\riii}{{\rm (iii)}}
\newcommand{\riv}{{\rm (iv)}}

\newcommand{\matriz}[1]{\ensuremath{\begin{bmatrix}#1\end{bmatrix}}}

\title[Complex solvmanifolds with trivial canonical bundle]{On the canonical bundle of complex solvmanifolds and applications to hypercomplex geometry}
\author{Adri\'an Andrada}
\thanks{Corresponding author: Adri\'an Andrada}
\email{adrian.andrada@unc.edu.ar}

\author{Alejandro Tolcachier}
\email{atolcachier@unc.edu.ar}

\address{FAMAF, Universidad Nacional de C\'ordoba and CIEM-CONICET, Av. Medina Allende s/n, X5000HUA C\'ordoba, Argentina}


\begin{document}

	\begin{abstract}
		We study complex solvmanifolds $\Gamma\backslash G$ with holomorphically trivial canonical bundle. We show that the trivializing section of this bundle can be either invariant or non-invariant by the action of $G$. First we characterize the existence of invariant trivializing sections in terms of the Koszul 1-form $\psi$ canonically associated to $(\g,J)$, where $\g$ is the Lie algebra of $G$, and we use this characterization to produce new examples of complex solvmanifolds with trivial canonical bundle. Moreover, we provide an algebraic obstruction, also in terms of $\psi$, for a complex solvmanifold to have trivial (or more generally holomorphically torsion) canonical bundle. Finally, we exhibit a compact hypercomplex solvmanifold $(M^{4n},\{J_1,J_2,J_3\})$ such that the canonical bundle of $(M,J_{\alpha})$ is trivial only for $\alpha=1$, so that $M$ is not an $\SL(n,\H)$-manifold.
	\end{abstract}
	
	\maketitle
	
	\section{Introduction}\label{S: intro}
	
	Given a complex manifold $(M,J)$ with $\dim_\C M=n$, its canonical bundle $K_M$ is defined as the $n$-th exterior power of its holomorphic cotangent bundle, and it is a holomorphic line bundle over $M$. This line bundle is holomorphically trivial when there exists a nowhere vanishing $(n,0)$-form which is holomorphic (or equivalently, closed). Complex manifolds with holomorphically trivial canonical bundle are important in differential geometry and other fields. For instance, compact K\"ahler manifolds $M$ with global Riemannian holonomy   contained in $\operatorname{SU}(n)$, have holomorphically trivial canonical bundle. More generally, any Calabi-Yau manifold (i.e., a compact K\"ahler manifold $M$ with $c_1(M)=0$ in $H^2(M,\R)$) has holomorphically torsion canonical bundle, that is, $K_M^{\otimes k}$ is trivial for some $k\in \N$. In theoretical physics, complex manifolds with holomorphically trivial canonical bundle appear in the study of the Hull-Strominger system. Indeed, in dimension 6, the solutions of this system occur in compact complex manifolds $M$ endowed with a special Hermitian metric (not necessarily K\"ahler) and with trivial $K_M$. According to \cite{To}, compact complex manifolds with holomorphically torsion canonical bundle have vanishing first Bott-Chern class, $c_1^{BC}=0$, and therefore they are examples of \textit{non-Kähler Calabi-Yau} manifolds.  
	
	A large family of compact complex manifolds with trivial canonical bundle is given by nilmanifolds $\Gamma\backslash G$ equipped with an invariant complex structure. Indeed, it was shown in \cite{BDV} that the simply connected nilpotent Lie group $G$ admits a nonzero left invariant holomorphic $(n,0)$-form $\sigma$ (with $\dim_\R G=2n$), by using a distinguished basis of left invariant $(1,0)$-forms provided by Salamon in \cite{Sal}. Since $\sigma$ is left invariant, it induces an invariant trivializing section of $K_{\Gamma\backslash G}$ for any lattice $\Gamma \subset G$.
	
	The next natural step is to study solvmanifolds $\Gamma\backslash G$ equipped with invariant complex structures (or complex solvmanifolds, for short). In this case, it is known that several different phenomena can occur. For instance:
	\begin{itemize}
		\item There are complex solvmanifolds which admit an invariant section of the canonical bundle, just as in the case of nilmanifolds. A classification of the Lie algebras associated to such solvmanifolds in dimension 6 is given in \cite{FOU}.
		\item There are $4$-dimensional complex solvmanifolds which do not have trivial canonical bundle. Indeed, Inoue surfaces are complex solvmanifolds (see for instance \cite{Hase}) with non trivial canonical bundle, since the only compact complex surfaces with trivial canonical bundle are complex tori, K3 surfaces and primary Kodaira surfaces.
		
	\end{itemize}
	
	\medskip 
	
	In this article we exhibit a different phenomenon concerning the canonical bundle of complex solvmanifolds. Indeed, in the next example we show that there exists a 4-dimensional complex solvmanifold $(\Gamma\backslash G,J)$ with trivial canonical bundle such that the trivializing section is not induced by a left invariant holomorphic $(2,0)$-form on $G$. This provides a counterexample to \cite[Proposition 2.1]{FOU}. 
	
	\begin{example}\label{ex: motivation}
		Let $H_3$ denote the $3$-dimensional Heisenberg group, which is considered as $\R^3$ equipped with the product
		\[ (x,y,z)\cdot (x',y',z')=(x+x',y+y',z+z'+\frac12 (xy'-yx')).\]
		Let us consider now the semidirect product $G=\R\ltimes_\varphi H_3$, where $\varphi: \R \to \operatorname{Aut}(H_3)\subset \operatorname{GL}(3,\R)$ is given by
		\[ \varphi(t)=\begin{bmatrix} \cos t & -\sin t & 0 \\ \sin t & \cos t & 0 \\ 0&0&1\end{bmatrix}. \]
		The Lie algebra $\g$ associated to $G$ has a basis $\{e_0,\ldots, e_3\}$ of left invariant vector fields such that 
		\[ [e_1,e_2]=e_3, \quad [e_0,e_1]=e_2, \quad [e_0,e_2]=-e_1;\]
		equivalently, the dual basis of $1$-forms $\{e^0,\ldots, e^3\}$ satisfies 
		\[  de^0=0, \quad de^1=e^{0}\wedge e^{2},\quad de^2=-e^{0}\wedge e^{1},\quad de^3=-e^{1}\wedge e^{2}. \]
		$G$ admits a left invariant complex structure $J$ given by $Je_0=e_3, Je_1=e_2$. A left invariant smooth section of the canonical bundle of $(G,J)$ is the $(2,0)$-form $\sigma=(e^0+ie^3) \wedge (e^1+ie^2)$. Thus, using the notation $e^{jk\cdots}=e^j \wedge e^k \wedge \cdots$, we have that $d\sigma=-i e^{023}-e^{013}\neq 0$. Therefore, $\sigma$ is not closed, and hence not holomorphic. However, if $t$ denotes the coordinate on the $\R$-factor we have that $dt=e^0$ and then the $(2,0)$-form $\tau:=\e^{it}\sigma$ is holomorphic, since $d\tau=0$.
		
		The Lie group $G$ admits the lattice $\Gamma=\{(2\pi k,m,n,\frac{p}{2})\mid k,m,n,p\in \Z\}$; note that $\Gamma$ is nilpotent and $\tau$ is $\Gamma$-invariant since $\e^{i(t+2\pi k)}=\e^{it}$. Therefore $\tau$  induces a nowhere vanishing closed $(2,0)$-form $\tilde{\tau}$ on the solvmanifold $(\Gamma\backslash G,J)$ and thus this solvmanifold has trivial canonical bundle. We point out that $(\Gamma\backslash G,J)$ is a primary Kodaira surface since it is biholomorphic to $(\Gamma\backslash (\R\times H_3), \tilde{J})$, where $\tilde{J}$ is induced by a left invariant complex structure on $\R\times H_3$, and noticing that $\Gamma$ is a lattice in both $G$ and $\R\times H_3$.
	\end{example}
	
	The previous example is the main motivation for this article. It shows that when studying the triviality of the canonical bundle of complex solvmanifolds we need to deal with the problem in two instances. First, given a complex solvmanifold $M=(\Gamma\backslash G,J)$ we have to determine if $M$ admits a trivializing section induced by a left invariant one on the Lie group so that $K_M$ is trivial. If this is not the case, then we must look for more general trivializing sections. Just as we did in Example \ref{ex: motivation}, we multiply a left invariant smooth section $\sigma$ of $K_{(G,J)}$ (which always exists) by a smooth function so as to get a holomorphic trivializing section of $K_{(G,J)}$. If this function is invariant by $\Gamma$ we obtain a trivializing section of $K_{(\Gamma\backslash G,J)}$, which is therefore is trivial. In this article we exhibit several examples of this situation. 
	
	Our main interest is to study complex solvmanifolds but whenever possible we prove results for compact quotients of general simply connected Lie groups equipped with left invariant complex structures.
	
	First, in \S \ref{S: inv}, given a (not necessarily solvable) Lie group $G$ endowed with a left invariant complex structure $J$, we tackle the problem of the existence of an invariant trivializing section of $K_{(G,J)}$. We show in Theorem \ref{theorem: inv} that such a section exists if and only if the Koszul 1-form $\psi$ on the Lie algebra $\g=\operatorname{Lie}(G)$ vanishes, where $\psi$ is defined
	by $\psi(x)=\Tr(J\ad x)-\Tr\ad (Jx)$, $x\in\g$. Applying this characterization we obtain new examples of complex solvmanifolds with trivial canonical bundle: for instance, when the complex structure is abelian (Corollary \ref{corollary: abelianFCT}). In addition, we show that if $(\Gamma\backslash G,J)$ admits an invariant trivializing section of some power of $K_{(\Gamma\backslash G,J)}$ then $K_{(\Gamma\backslash G,J)}$ itself admits an invariant trivializing section. 
	
	In \S \ref{S: groups} we prove first that two trivializing sections of the canonical bundle of a Lie group $G$ equipped with a left invariant complex structure differ by a nowhere vanishing holomorphic function on $G$ (Lemma \ref{lemma: holomorfas}). This implies that the trivializing sections of the canonical bundle of a compact complex manifold $\Gamma\backslash G$ are either all invariant or all non invariant (Corollary \ref{corollary: cpt uniqueness}). Then we prove that any simply connected solvable Lie group $G$ equipped with a left invariant complex structure $J$ has trivial canonical bundle (Theorem \ref{theorem:non-invariant}). 
	
	In \S \ref{S: solvmflds}, our goal is to give new examples of complex solvmanifolds with trivial canonical bundle, in cases where there are no invariant sections. We first provide an algebraic obstruction in terms of the Koszul 1-form $\psi$, which holds for any Lie group $G$ with a left invariant complex structure. Indeed, 
	we show in Theorem \ref{theorem: obstruction} that if a compact complex manifold $\Gamma\backslash G$ has trivial (or more generally holomorphically torsion) canonical bundle then $\psi$ vanishes on the commutator ideal $[\g,\g]$ where $\g=\operatorname{Lie}(G)$. We use this condition to reobtain the known fact that compact semisimple Lie groups with a left invariant complex structure do not have holomorphically torsion canonical bundle (Proposition \ref{proposition:compactos}). This obstruction also provides us with a helpful insight to find an explicit trivializing section of the canonical bundle of some complex solvmanifolds (Proposition \ref{proposition: tau explicit}). We also apply this construction in order to exhibit some new examples.
	
	In the last section we consider a Lie group $G$ equipped with a left invariant hypercomplex structure $\{J_1,J_2,J_3\}$ and we study the triviality of the canonical bundle of the complex manifolds $(G,J_\alpha)$, $\alpha=1,2,3$.
	First  we prove in Theorem \ref{theorem: hcpx} that if $\{J_1,J_2,J_3\}$ is a left invariant hypercomplex structure on $G$ and if $(G,J_\alpha)$ admits a left invariant trivializing section of its canonical bundle for some $\alpha=1,2,3$, then the canonical bundle of $(G,J_\beta)$ is trivial for all $\beta=1,2,3$, and the same happens for any hypercomplex compact quotient $\Gamma\backslash G$. Next we show that this does not necessarily hold for hypercomplex solvmanifolds if the trivializing section of $(\Gamma\backslash G,J_\alpha)$ is not invariant. Indeed, in Example \ref{ex: hcpx} we exhibit 8-dimensional hypercomplex solvmanifolds $(\Gamma\backslash G,\{J_1,J_2,J_3\})$ such that $(\Gamma\backslash G,J_1)$ has trivial canonical bundle but the canonical bundles of $(\Gamma\backslash G,J_2)$ and $(\Gamma\backslash G,J_3)$ are both non trivial. These examples are not $\SL(2,\mathbb{H})$-manifolds hence they provide a negative answer to a question by Verbitsky in \cite{Ver}.

	\smallskip
	
	\textbf{Acknowledgments.} 
	The authors are grateful to Daniele Angella for his useful comments and suggestions and also to the anonymous reviewers for their careful reading of the manuscript.
	This work was partially supported by CONICET, SECyT-UNC and FONCyT (Argentina).

	\section{Preliminaries}\label{S: prelim}
	
	An almost complex structure on a differentiable manifold $M$ is an automorphism
	$J$ of the tangent bundle $TM$ satisfying $J^2=-\I $, where $\I$ is the identity endomorphism of $TM$. Note that the existence of an almost complex structure on $M$ forces the dimension of $M$ to be even, say $\dim_\R M=2n$. The almost complex structure $J$ is called \textit{integrable} when it satisfies the condition $N_{J}\equiv 0$, where $N_J$ is the Nijenhuis tensor given by:
	\begin{equation}\label{eq:nijenhuis}
		N_{J}(X,Y) =
		[X,Y]+J([JX,Y]+[X,JY])-[JX,JY],
	\end{equation}
	for $X,Y$ vector fields on $M$. An integrable almost complex structure is called simply a complex structure on $M$. According to the well-known Newlander-Nirenberg theorem, a complex structure on $M$ is equivalent to the existence of a holomorphic atlas on $M$, so that $(M,J)$ can be considered as a complex manifold of complex dimension $n$. 
	
	If $(M,J)$ is a complex manifold with $\dim_\C M=n$ its canonical bundle is defined as 
	\[ K_M=\alt^n {\mathcal T}^*_M,  \]
	where ${\mathcal T}^*_M$ is the holomorphic cotangent bundle of $M$. This is a holomorphic line bundle on $M$, and it is holomorphically trivial if and only if there exists a nowhere vanishing holomorphic $(n,0)$-form defined on $M$. In this article, by trivial canonical bundle we will always mean holomorphically trivial canonical bundle. 
	
	Note that if $\sigma$ is a $(n,0)$-form on $M$ then $\sigma$ is holomorphic if and only if it is closed, since $d\sigma=\partial \sigma+\overline{\partial}\sigma$ and $\partial \sigma$ is a $(n+1, 0)$-form, thus $\partial \sigma=0$. 
	
	We observe first that the existence of a trivializing section of the canonical bundle has some topological consequences on a compact complex manifold.
	
	\begin{proposition}\label{proposition:betti-n}
		Let $(M,J)$ be a compact complex manifold with trivial canonical bundle and $\dim_\R M=2n$. Then the $n$-th Betti number $b_n(M)$ satisfies $b_n(M)\geq 2$.  
	\end{proposition}
	
	
	\begin{proof}
		We follow the lines of \cite[Proposition 2.5]{FOU}. Let $\tau$ be a nowhere vanishing holomorphic $(n,0)$-form on $M$, therefore $\tau\wedge \bar{\tau}$ is a nonzero multiple of a real volume form on $M$. Let us decompose it as $\tau=\tau_1+i\tau_2$. Since $\tau$ is closed, we have that $d\tau_1=0=d\tau_2$. Therefore, they define de Rham cohomology classes $[\tau_1],[\tau_2]\in H^n_{dR}(M,\R)$. These two classes are linearly independent. Indeed, if we assume otherwise then there exist $a,b\in\R$ with $a^2+b^2\neq 0$ such that $a\tau_1+b\tau_2=d\eta$ for some $(n-1)$-form $\eta$. We have two cases, according to the parity of $n$. 
		
		$\ri$ Case $n$ odd: in this case we have \[0\neq \tau\wedge \bar{\tau}=(\tau_1\wedge\tau_1 + \tau_2\wedge\tau_2)+i(-\tau_1\wedge\tau_2+\tau_2\wedge\tau_1)=-2i (\tau_1 \wedge \tau_2).\]
		We compute next 
		\[ d(\eta\wedge (-b\tau_1+a\tau_2))=(a\tau_1+b\tau_2)\wedge (-b\tau_1+a\tau_2)=(a^2+b^2)(\tau_1\wedge\tau_2).\]
		Integrating over $M$ we obtain, due to Stokes' theorem, $\displaystyle{0=(a^2+b^2)\int_M \tau_1\wedge\tau_2}$, which is a contradiction.
		
		$\rii$ Case $n$ even: in this case we have
		\[ 0\neq \tau\wedge \bar{\tau}=(\tau_1\wedge\tau_1 + \tau_2\wedge\tau_2)+i(-\tau_1\wedge\tau_2+\tau_2\wedge\tau_1)=\tau_1 \wedge\tau_1 + \tau_2 \wedge \tau_2.\]
		It follows from $0=\tau\wedge \tau$ that $\tau_1\wedge\tau_1=\tau_2 \wedge \tau_2$ and $\tau_1\wedge\tau_2=0$. In particular, $0\neq \tau\wedge\bar{\tau}=2\tau_1\wedge \tau_1$. We compute next
		\[ d(\eta\wedge (a\tau_1 +b\tau_2))=(a\tau_1+b\tau_2)\wedge (a\tau_1+b\tau_2)=(a^2+b^2)\tau_1\wedge\tau_1.\]
		Again, integrating over $M$ we obtain a contradiction.
		
		Therefore we obtain that $b_n(M)\geq2$.
	\end{proof}
	
	\smallskip 
	
	Other important holomorphic line bundles over the complex manifold $(M,J)$
	are given by the tensor powers of the canonical bundle:
	\[ K_M^{\otimes k}=K_M\otimes \cdots \otimes K_M \quad \text{($k$ times)}.\]
	Following \cite{To}, we will say that a complex manifold $(M,J)$ is \textit{holomorphically torsion} if $K_M^{\otimes k}$ is holomorphically trivial for some $k\geq 1$. The triviality of this holomorphic bundle can be understood as follows. 
	
	For any complex manifold $M$ the Dolbeault operator $\bar{\partial}$ can be extended to a differential operator $\bar{\partial}_k:\Gamma(K_M^{\otimes k})\to \Gamma((T^*M)^{0,1}\otimes K_M^{\otimes k}$), where $\Gamma(\cdot)$ denotes the space of smooth sections. Indeed, since $\alt^{n,1}(M)\cong (T^* M)^{0,1} \otimes K_M$ we define recursively: $\bar{\partial}_1=\bar{\partial}$ and for $k\geq 2$,
	\[ \bar{\partial}_k(\sigma\otimes s)=\bar{\partial}\sigma\otimes s+\sigma\otimes \bar{\partial}_{k-1} s,\] where $\sigma\in \Gamma(K_M)$ and $s\in \Gamma(K_M^{\otimes k-1})$. This differential operator satisfies the Leibniz rule $\bar{\partial}_k(fs)=\bar{\partial}f\otimes s+f \bar{\partial}_k s$ for any $f\in C^{\infty}(M,\C)$ and $s\in\Gamma(K_M^{\otimes k})$.
	
	The holomorphic bundle $K_M^{\otimes k}$ is trivial if and only if there exists a nowhere vanishing section $s\in\Gamma(K_M^{\otimes k})$ such that $\bar{\partial}_k s=0$.
	
	\medskip
	
	A Hermitian structure on a smooth manifold $M$ is a pair $(J,g)$ of a complex structure $J$ and a Riemannian metric $g$ compatible with $J$, that is, $g (JX, JY) =g( X, Y) $ for all vector fields $X, Y$ on $M$, or equivalently, $g (JX,Y)=-g(X,JY)$. 
	
	\smallskip
	
	\subsection{Solvmanifolds}
	
	A discrete subgroup $\Gamma$ of a Lie group $G$ is called a \textit{lattice} if the quotient $\Gamma\backslash G$ has finite volume. According to \cite{Mil}, if such a lattice exists then the Lie group must be unimodular, that is, it carries a bi-invariant Haar measure. This is equivalent, when $G$ is connected, to $\Tr( \ad x)=0$ for all $x\in \g=\operatorname{Lie}(G)$ (in this case, $\g$ is called unimodular as well). When $\Gamma\backslash G$ is compact the lattice $\Gamma$ is said to be uniform. It is well known that when $G$ is solvable then any lattice is uniform \cite[Theorem 3.1]{Rag}. 
	
	Assume that $G$ is simply connected and $\Gamma$ is a uniform lattice in $G$. 
	The quotient $\Gamma\backslash G$ is called a solvmanifold if $G$ is solvable and a nilmanifold if $G$ is nilpotent, and it follows that $\pi_1(\Gamma\backslash G)\cong \Gamma$. 
	Furthermore, the diffeomorphism class of solvmanifolds is determined by the isomorphism class of the corresponding lattices, as the following result shows:
	
	\begin{theorem}\cite{Mo}\label{solv-isom}
		If $\Gamma_1$ and $\Gamma_2$ are lattices in simply connected solvable Lie groups 
		$G_1$ and $G_2$, respectively, and $\Gamma_1$ is isomorphic to $\Gamma_2$, then $\Gamma_1 \backslash G_1$ is diffeomorphic to $\Gamma_2 \backslash G_2$.
	\end{theorem}
	
	
	
	Note that in any fixed dimension only countably many non-isomorphic simply connected Lie groups admit lattices, according to \cite{Milo} (for the solvable case) and \cite{Wi} (for the general case).
	
	\medskip
	
	Let $G$ be a simply connected solvable Lie group, and $N$ the nilradical of $G$ (i.e., the connected closed Lie subgroup of $G$ whose Lie algebra is the nilradical $\n$ of $\g$). Moreover, $[G,G]$ is the connected closed Lie subgroup with Lie algebra $[\g,\g]$. As $G$ is solvable, $[G,G]\subset N$ so $G/N$ is abelian, and from the long exact sequence of homotopy groups associated to the fibration $N\to G\to G/ N$ it follows that $G/N$ is simply connected. Therefore $G/N \cong \R^k$ for some $k\in\N$ and $G$ satisfies the short exact sequence \[1\to N\to G\to \R^k\to 1.\] $G$ is called \textit{splittable} if this sequence splits, that is, there is a right inverse homomorphism of the projection $G\to \R^k$. This condition is equivalent to the existence of a homomorphism $\phi:\R^k\to \Aut(N)$ such that $G$ is isomorphic to the semidirect product $\R^k \ltimes_\phi N$.
	
	Following \cite{Yamada}, a lattice $\Gamma$ of a splittable solvable Lie group $\R^k\ltimes_{\phi} N$ will be called \textit{splittable} if it can be written as $\Gamma=\Gamma_1\ltimes_{\phi} \Gamma_2$ where $\Gamma_1\subset\R^k$ and $\Gamma_2\subset N$ are lattices of $\R^k$ and $N$ respectively. Also in \cite{Yamada} there is a criterion to determine the existence of splittable lattices in splittable solvable simply connected Lie groups.
	
	\begin{theorem}\cite{Yamada}\label{theorem: yamada}
		Let $G=\R^k \ltimes_\phi N$ be a simply connected splittable solvable Lie group, where $N$ is the nilradical of $G$. If there exist a rational basis $\mathcal{B}=\{X_1,\ldots,X_n\}$ of $\n$ and a basis $\{t_1,\ldots,t_k\}$ of $\R^k$ such that the coordinate matrix of $d(\phi(t_j))_{1_N}$ in the basis $\mathcal{B}$ is an integer unimodular matrix for all $1\leq j\leq k$ then $G$ has a splittable lattice of the form $\Gamma=\text{span}_\Z\{t_1,\ldots,t_k\}\ltimes_\phi \exp^N(\text{span}_\Z \{X_1,\ldots,X_n\})$.
	\end{theorem} 
	
	When $k=1$ the simply connected solvable splittable Lie group $G=\R \ltimes_{\phi} N$ is called \textit{almost nilpotent}. In this case, every lattice is splittable due to \cite{Bock}. If $N$ is abelian, i.e. $N=\R^n$, then $G$ (and its corresponding Lie algebra) is called \textit{almost abelian}. 
	
	In the examples in the forthcoming sections, we will begin with a Lie algebra $\g=\R^k\ltimes_{\varphi}\n$. In order to apply Theorem \ref{theorem: yamada} we need to determine the associated simply connected Lie group $G$. Let $N$ denote the simply connected nilpotent Lie group with Lie algebra $\n$. Since $\exp:\n\to N$ is a diffeomorphism, we may assume that the underlying manifold of $N$ is $\n$ itself with the group law $x\cdot y=Z(x,y)$, where $Z(x,y)$ is the polynomial map given by the Baker-Campbell-Hausdorff formula: $\exp(x)\exp(y)=\exp(Z(x,y))$. Therefore, with this assumption, we have that $\exp:\n\to N$ is simply the identity map on $\n$ and moreover, $\Aut(\n)=\Aut(N)$.
	
	Let $\{t_1,\ldots,t_k\}$ be a basis of $\R^k$ and denote $B_j=\varphi(t_j)\in \Der(\n)$. Then, $\exp(B_j)\in \Aut(N)$ and using \cite[Theorem 4.2]{Bock} we have that $G=\R^k\ltimes_{\phi} N$, where $\phi:\R^k\to \Aut(N)$ is the Lie group homomorphism given by \[ \phi\left(\sum_{j=1}^k  x_jt_j\right)=\exp(x_1 B_1+\cdots+x_k B_k)=\exp(x_1 B_1)\exp(x_2 B_2)\cdots \exp(x_k B_k).\] Here  $\exp$ denotes the matrix exponential after identification of $\n\cong \R^{\dim \n}$ choosing a basis of $\n$. 
	
	Note that, in the notation of Theorem \ref{theorem: yamada}, we have that  $d(\phi(t_j))_{1_N}=\exp(B_j)=\exp(\varphi(t_j))$. Hence, in order to find lattices we need a basis $\{t_1,\ldots,t_k\}$ such that
	$[\exp(\varphi(t_j))]_\B$ is an integer unimodular matrix in the rational basis $\B$ of $\n$, for all $1\leq j\leq k$.

	\smallskip
	
	We move on now to consider invariant geometric structures on solvmanifolds. 
	
	Let $G$ be a connected Lie group with Lie algebra $\g$. A complex structure $J$ on $G$ is said to be left invariant if left translations by elements of $G$ are holomorphic maps. In this case $J$ is determined by the value at the identity of $G$. Thus, a left invariant complex structure on $G$ amounts to a complex structure on its Lie algebra $\g$, that is, a real linear transformation $J$ of $\g$ satisfying $J^2 = -\I$ and $N_J(x, y)=0$ for all $x, y$ in $\g$, with $N_J$ defined as in \eqref{eq:nijenhuis}. A Riemannian metric $g$ on $G$ is called left invariant when left translations are isometries. Such a metric $g$ is determined by its value $g_e=\pint$ at the identity $e$ of $G$, that is, $\pint$ is a positive definite inner product on $T_e G=\g$. A Hermitian structure $(J,g)$ on $G$ is left invariant when both $J$ and $g$ are left invariant. The corresponding pair $(J,\pint)$ is called a Hermitian structure on $\g$.
	
	We observe that left invariant geometric structures defined on $G$ induce naturally geometric structures on $\Gamma\backslash G$, with $\Gamma$ a lattice in $G$, which are called \textit{invariant}. For instance, a left invariant complex structure (respectively, Riemannian metric) on $G$ induces a unique complex structure (respectively, Riemannnian metric) on $\Gamma\backslash G$ such that the canonical projection $G\to \Gamma\backslash G$ is a local biholomorphism (respectively, local isometry). In this article, a solvmanifold equipped with an invariant complex structure will be called simply a \textit{complex solvmanifold}.
	
	\medskip 
	
	\section{Complex solvmanifolds with trivial canonical bundle via invariant sections} \label{S: inv}
	In this section we deal with the existence of nowhere vanishing left invariant closed $(n,0)$-forms on $2n$-dimensional Lie groups equipped with a left invariant complex structure, equivalently we study the existence of nonzero closed $(n,0)$-forms on the corresponding Lie algebras.
	
	First we characterize the existence of such a form in algebraic terms. In order to do so, we need the following notion which will play a crucial role throughout the article:
	
	\begin{definition}
		Given a Lie algebra $\g$ equipped with a complex structure $J$, the \textit{Koszul $1$-form} $\psi\in\g^*$ is defined as:
		\begin{equation}\label{eq: psi} 
			\psi(x)=\Tr (J\ad x)-\Tr \ad (Jx), \quad x\in\g.
		\end{equation}
		This form was introduced by Koszul in \cite{Koszul} (see also \cite{Grant}) in the context of $G$-invariant complex structures on homogeneous spaces $G/H$ (compare with $\theta^1$ in \cite[Proposition 4.1]{V}).
	\end{definition}

	

	\begin{theorem}\label{theorem: inv}
		Let $\g$ be a $2n$-dimensional Lie algebra with an almost complex structure $J$. Let $\sigma\in \alt^{n,0} \g^*$ be a nonzero $(n,0)$-form on $\g$. Then $d\sigma=0$ if and only if $J$ is integrable and $\psi \equiv 0$. 
	\end{theorem}
	
	\begin{proof}
		Let $\{u_1,\ldots,u_n,v_1,\ldots,v_n\}$ be a $J$-adapted basis of $\g$, that is, $Ju_j=v_j$ for all $j$.
		Since $\dim_\C \alt^{n,0} \g^*=1$, we may assume that $\sigma=(u^1+iv^1)\wedge\cdots\wedge(u^n+iv^n)$.
		
		The Lie brackets of $\g$ can be written in terms of the basis above by  
		\begin{align*}
			[u_j,u_k]=\sum_{\l=1}^n a_{jk}^\l u_\l+\sum_{\l=1}^n b_{jk}^\l v_\l, \quad 
			[u_j,v_k]=\sum_{\l=1}^n c_{jk}^\l u_\l+\sum_{\l=1}^n d_{jk}^\l v_\l, \quad
			[v_j,v_k]=\sum_{\l=1}^n e_{jk}^\l u_\l+\sum_{\l=1}^n f_{jk}^\l v_\l,
		\end{align*}
		with $a_{kj}^{\l}=-a_{jk}^{\l}$, $b_{kj}^{\l}=-b_{jk}^{\l}$, $e_{kj}^{\l}=-e_{jk}^{\l}$ and $f_{kj}^{\l}=-f_{jk}^{\l}$. 
		Accordingly, for $1\leq \l\leq n$ we have
		\begin{align*} 
			du^\l&=-\sum_{j,k=1}^n \left(\tfrac12 a_{jk}^\l\,  u^{jk}+ c_{jk}^\l\,u^j\wedge v^k+\tfrac12 e_{jk}^\l\, v^{jk}\right),\\
			dv^\l&=-\sum_{j,k=1}^n \left(\tfrac12 b_{jk}^\l\, u^{jk}+d_{jk}^\l\, u^j\wedge v^k+\tfrac12 f_{jk}^\l\, v^{jk}\right).
		\end{align*}
		
		Let us set $\gamma_j:=u^j+iv^j$ for all $j$, so that $\sigma=\gamma_1\wedge\cdots\wedge \gamma_n$. Next we compute $d\gamma_\l$ in terms of $\gamma_j$ and $\bar{\gamma_j}$. First we note that $2u^j=(\gamma_j+\bar{\gamma_j})$ and $2v^j=-i(\gamma_j-\bar{\gamma_j})$ imply that
		\begin{align*}
			4u^{jk}=(\gamma_j+\bar{\gamma_j})\wedge (\gamma_k+\bar{\gamma_k})&=\gamma_{jk}+\gamma_{\bar{j}k}+\gamma_{j\bar{k}}+\gamma_{\bar{j}\,\bar{k}}\\
			4u^j\wedge v^k=-i(\gamma_j+\bar{\gamma_j})\wedge (\gamma_k-\bar{\gamma_k})&=-i(\gamma_{jk}+\gamma_{\bar{j}k}-\gamma_{j\bar{k}}-\gamma_{\bar{j}\,\bar{k}})\\
			4v^{jk}=-(\gamma_j-\bar{\gamma_j})\wedge (\gamma_k-\bar{\gamma_k})&=-(\gamma_{jk}-\gamma_{\bar{j}k}-\gamma_{j\bar{k}}+\gamma_{\bar{j}\,\bar{k}})
		\end{align*}
		Using these identities it follows that 
		\begin{align}\label{eq: dgamma_l}
			d\gamma_\l&=-\sum_{j,k=1}^n  \tfrac12 (a_{jk}^\l+ib_{jk}^\l) u^{jk}+(c_{jk}^\l+id_{jk}^\l) u^j\wedge v^k+\tfrac12 (e_{jk}^\l+if_{jk}^\l) v^{jk}\\
			&=-\tfrac14 \sum_{j,k=1}^n\Bigg( \left(\tfrac12 a_{jk}^\l+d_{jk}^\l-\tfrac12 e_{jk}^\l+i\left(\tfrac12 b_{jk}^\l-c_{jk}^\l-\tfrac12 f_{jk}^\l\right)\right) \gamma_{jk}\nonumber\\
			&\quad+\left(\tfrac12 a_{jk}^\l+d_{jk}^\l+\tfrac12 e_{jk}^\l+i\left(\tfrac12 b_{jk}^\l-c_{jk}^\l+\tfrac12 f_{jk}^\l\right)\right) \gamma_{\bar{j}k}\nonumber\\
			&\quad+\left(\tfrac12 a_{jk}^\l-d_{jk}^\l+\tfrac12 e_{jk}^\l+i\left(\tfrac12 b_{jk}^\l+c_{jk}^\l+\tfrac12 f_{jk}^\l\right)\right)\gamma_{j\bar{k}} \nonumber \\
			&\quad+\left(\tfrac12 a_{jk}^\l-d_{jk}^\l-\tfrac12 e_{jk}^\l+i\left(\tfrac12 b_{jk}^\l+c_{jk}^\l-\tfrac12 f_{jk}^\l\right)\right)\gamma_{\bar{j}\,\bar{k}}\Bigg) \nonumber 
		\end{align}
		We use the previous expression in order to compute $d\sigma$, where we use the unconventional but shorter notation $\dfrac{\sigma}{\gamma_l}= \gamma_1\wedge \cdots  \wedge \gamma_{\l-1}  \wedge\gamma_{\l+1} \wedge \cdots \wedge \gamma_n$: 
		\begin{align*}
			d\sigma&=\sum_{\l=1}^{n} (-1)^{\l+1} \, \gamma_1\wedge \cdots  \wedge \gamma_{\l-1} \wedge d\gamma_\l \wedge\gamma_{\l+1} \wedge \cdots \wedge \gamma_n\\
			&=-\tfrac 14\sum_{j,\l} \left(\tfrac12a_{j\l}^\l+d_{j\l}^\l+\tfrac12 e_{j\l}^\l+i\left(\tfrac12 b_{j\l}^\l-c_{j\l}^\l+\tfrac12 f_{j\l}^\l\right)\right) \bar{\gamma_j} \wedge \sigma \  \\
			&\quad + \tfrac 14\sum_{k,\l} \left(\tfrac12 a_{\l k}^\l-d_{\l k}^\l+\tfrac12 e_{\l k}^\l+i\left(\tfrac12 b_{\l k}^\l+c_{\l k}^\l+\tfrac12 f_{\l k}^\l\right)\right) \bar{\gamma_k} \wedge \sigma \\\
			&\quad +\tfrac 14 \sum_{j,k,\l=1}^n (-1)^\l  \left(\tfrac12 a_{jk}^\l-d_{jk}^\l-\tfrac12 e_{jk}^\l+i\left(\tfrac12 b_{jk}^\l+c_{jk}^\l-\tfrac12 f_{jk}^\l\right)\right)\bar{\gamma_j} \wedge\bar{\gamma_k} \wedge \frac{\sigma}{\gamma_\l} \\\
			&=\tfrac14 \sum_{j,\l} (a_{\l j}^\l-(d_{j\l}^\l+d_{\l j}^\l)+e_{\l j}^\l) +i(b_{\l j}^\l+(c_{\l j}^\l+c_{j\l}^\l)+f_{\l j}^\l)) \bar{\gamma_j} \wedge \sigma \\\
			\quad&+\tfrac14\sum_{j<k}\sum_{\l} ((a_{jk}^\l+(-d_{jk}^\l+d_{kj}^\l)-e_{jk}^\l)+i(b_{jk}^\l+(c_{jk}^\l-c_{kj}^\l)-f_{jk}^\l)) \bar{\gamma_j}\wedge\bar{\gamma_k}\wedge\frac{\sigma}{\gamma_\l}. 
		\end{align*} 
		Since $\{\bar{\gamma_j}\wedge\sigma\mid 1\leq j\leq n\}\cup \{\bar{\gamma_j}\wedge\bar{\gamma_k}\wedge\frac{\sigma}{\gamma_\l}\mid j<k, 1\leq \l\leq n\}$ is linearly independent, we see that $d\sigma=0$ if and only if
		\begin{align}\label{eq: 1° dsigma=0}
			\sum_{\l=1}^n a_{\l j}^\l-d_{j\l}^\l-d_{\l j}^\l+e_{\l j}^\l=0,& \quad \sum_{\l=1}^n b_{\l j}^\l+c_{\l j}^\l+c_{j\l}^\l+f_{\l j}^\l=0,
			\quad 1\leq j\leq n,\\ \label{eq: 2° dsigma=0}
			e_{jk}^\l=a_{jk}^\l-d_{jk}^\l+d_{kj}^\l&,\; f_{jk}^\l=b_{jk}^\l+c_{jk}^\l-c_{kj}^\l,\quad j<k,\; 1\leq \l \leq n. 
		\end{align}
		On the other hand, it is well known that the integrability of $J$ is equivalent to 
		\[ d(\alt^{1,0}\g^\ast_\C) \subseteq \alt^{2,0}\g^\ast_\C\oplus \alt^{1,1}\g^\ast_\C, \] where $\g_\C$ denotes the complexification of $\g$ and the bidegrees are induced by $J$.
		
		Therefore, $J$ is integrable if and only if the coefficient of $\bar{\gamma_j}\wedge\bar{\gamma_k}$ in $d\gamma_{\l}$ vanishes for all $j,k,\l$. It follows from \eqref{eq: dgamma_l} that this happens if and only if
		\begin{align*}
			e_{jk}^\l=a_{jk}^\l-d_{jk}^\l+d_{kj}^\l, \quad f_{jk}^\l=b_{jk}^\l+c_{jk}^\l-c_{kj}^\l,\quad j<k, \; 1\leq \l \leq n,
		\end{align*}
		which is exactly \eqref{eq: 2° dsigma=0}.
		
		Next, using the inner product $\la\cdot,\cdot\ra$ on $\g$ defined by decreeing the basis $\{u_1,\ldots,u_n,v_1,\ldots,v_n\}$ orthonormal we have that 
		\begin{align*}
			\Tr(J \ad u_j)=\sum_{\l=1}^n (b_{\l j}^\l+c_{j \l}^\l)&, \quad \Tr(J \ad v_j)=\sum_{\l=1}^n (d_{\l j}^\l-e_{\l j}^\l),\\
			-\Tr \ad v_j =\sum_{\l=1}^n (c_{\l j}^\l+f_{\l j}^\l)&, \quad \Tr \ad u_j=\sum_{\l=1}^n (a_{j\l}^\l+d_{j\l}^\l).\end{align*}
		Hence, \eqref{eq: 1° dsigma=0} can be written as, for $1\leq j\leq n$,
		\[ -\Tr \ad u_j-\Tr(J\ad v_j)=0, \qquad 
		-\Tr \ad v_j+\Tr(J\ad u_j)=0. \]
		Then \eqref{eq: 1° dsigma=0} is equivalent to $\psi(u_j)=\psi(v_j)=0$. Thus, $d\sigma=0$ if and only if $J$ is integrable and $\psi\equiv0$.
	\end{proof}
	
	\smallskip 
	
	In the unimodular case we obtain the following characterization.
	
	\begin{corollary}\label{corollary: inv unimod}
		Let $\g$ be a $2n$-dimensional unimodular Lie algebra with a complex structure $J$. Then $\g$ admits a closed non-vanishing $(n,0)$-form   if and only if $\Tr(J\ad x)=0$ for all $x\in\g$. 
	\end{corollary}
	
	\smallskip
	
	\begin{remark}\label{remark: dsigma}
		It follows from the proof of Theorem \ref{theorem: inv} that if $J$ is integrable then\footnote{Cf. \cite[Lemma 3]{GZ}} 
		\begin{align*}
			d\sigma&= \tfrac14 \sum_{j,\l} (a_{\l j}^\l-d_{j\l}^\l-d_{\l j}^\l+e_{\l j}^\l)+i(b_{\l j}^\l+c_{\l j}^\l+c_{j\l}^\l+f_{\l j}^\l) \; \bar{\gamma_j} \wedge \sigma \\
			&=\tfrac14 \sum_j \left(-\Tr(\ad u_j)-\Tr(J\ad v_j))+i(\Tr(J \ad u_j)-\Tr(\ad v_j))\right)\; \bar{\gamma_j} \wedge \sigma\\
			&=\tfrac14 \sum_j (-\psi(v_j)+i \psi(u_j)) \;\bar{\gamma_j} \wedge \sigma. 
		\end{align*} 
	\end{remark}
	
	\ 
	
	When $\g$ is unimodular and $J$ is integrable, the vanishing of the Koszul 1-form $\psi$ can also be understood in terms of the complexification $\g_\C$ of $\g$ as the  following proposition shows. Recall that $\g_\C=\g^{1,0}\oplus\g^{0,1}$, where $\g^{1,0}$ (respectively, $\g^{0,1}$) is the $i$-eigenspace (respectively, $(-i)$-eigenspace) of the $\C$-linear extension $J^\C:\g_\C\to\g_\C$, and they are given by 
	\[
	\g^{1,0}=\{x-iJx\mid x\in \g\},\quad 
	\g^{0,1}=\{x+iJx\mid x\in \g\}.
	\] Both $\g^{1,0}$ and $\g^{0,1}$ are Lie subalgebras of $\g_\C$ due to the integrability of $J$.
	
	\begin{proposition}
		Let $(\g,J)$ be a $2n$-dimensional unimodular Lie algebra equipped with a complex structure. Then $(\g,J)$ has a nonzero closed $(n,0)$-form if and only if $\g^{1,0}$ (or $\g^{0,1}$) is unimodular.
	\end{proposition}
	
	\begin{proof}
		We will show the equivalence only for $\g^{1,0}$. The computations for $\g^{0,1}$ are completely analogous. Given a Hermitian inner product $\pint$ on $\g$, it can be extended to a complex inner product on $\g_\C$ satisfying $\la a+ib,c+id\ra=\la a,c\ra +\la b,d\ra-i(\la a,d\ra-\la b,c\ra)$. Thus, if $\{e_j\}_{j=1}^{2n}$ is an orthonormal basis of $\g$ such that $Je_{2j-1}=e_{2j}$ then $\{\frac{1}{\sqrt{2}}(e_{2j-1}- i e_{2j})\}_{j=1}^n$ is an orthonormal basis of $\g^{1,0}$.
		
		Now, consider $x-iJx\in\g^{1,0}$. We can decompose $\ad(x-iJx)$ with respect to the decomposition $\g_\C=\g^{1,0}\oplus \g^{0,1}$ as
		\[ \ad(x-iJx)=\begin{array}{|c|c|} A_x & * \\ \hline 0 & B_x \end{array}. \]
		Next we compute
		\begin{align*}
			\Tr A_x&=\tfrac12 \sum_{j=1}^n \la [x-iJx, e_{2j-1}-i e_{2j}], e_{2j-1}-i e_{2j}\ra\\
			&=\frac12 \sum_{j=1}^n \la [x,e_{2j-1}]-[Jx,e_{2j}]-i([x,e_{2j}]+[Jx,e_{2j-1}]), e_{2j-1}-i e_{2j}\ra \\
			&=\tfrac12 \left(\Tr \ad x-i\Tr \ad (Jx)-\Tr(J \ad(Jx))-i \Tr(J \ad x) \right)\\
			&=\tfrac12 (\Tr \ad x-\Tr(J \ad(Jx))-\tfrac{i}{2}(\Tr \ad(Jx)+\Tr(J\ad x)). 
		\end{align*}
		Therefore, $\Tr \ad(x-iJx)=0$ on $\g^{1,0}$ if and only if $\Tr (J\ad x)=-\Tr \ad (Jx)$ on $\g$. In particular, as $\g$ is unimodular, it follows that $\g^{1,0}$ is unimodular if and only if $\Tr(J \ad x)=0$, and the statement follows from Corollary \ref{corollary: inv unimod}.
	\end{proof}
	

	\smallskip

	There are many known examples of complex compact quotients $\Gamma\backslash G$ of a simply connected Lie group $G$ by a discrete subgroup $\Gamma$ admitting an invariant trivializing section of the canonical bundle. For instance:  
	\begin{itemize}
		\item when $G$ is a complex Lie group, that is, the complex structure $J:\g\to \g$ satisfies $J\ad x=\ad(Jx)$ for all $x\in \g$ (such a complex structure is called \textit{bi-invariant}) since $(G,J)$ is complex parallelizable,
		\item when $G$ is nilpotent \cite{BDV},
		\item some complex almost abelian solvmanifolds \cite{FP}.   
	\end{itemize}
	
	All these examples are easily recovered  
	using Theorem \ref{theorem: inv}.
	
	\medskip
	
	We consider next a special family of left invariant complex structures on Lie groups. If an almost complex structure $J$ on a Lie algebra $\g$ satisfies $[Jx,Jy]=[x,y]$ for all $x,y\in\g$ then it is immediate to verify that $J$ is integrable. Such a complex structure is called \textit{abelian}. They were introduced in \cite{BDM} and they have proved very useful in different contexts in differential and complex geometry. Abelian complex structures can only  occur on $2$-step solvable Lie algebras (see for instance \cite{ABDO}).  
	
	In the next result we show the existence of a left invariant trivializing section of the canonical bundle of a unimodular Lie group equipped with an abelian complex structure. As usual we state the result at the Lie algebra level.
	
	\begin{corollary}\label{corollary: abelianFCT}
		A $2n$-dimensional Lie algebra $\g$ equipped with an abelian complex structure $J$ has a nonzero closed $(n,0)$-form if and only if $\g$ is unimodular. In particular, any complex solvmanifold equipped with an abelian complex structure has trivial canonical bundle.
	\end{corollary} 
	
	\begin{proof}
		The fact that $J$ is abelian is equivalent to $[x,Jy]=-[Jx,y]$ for all $x,y\in\g$. Hence, $\ad(x) J=-\ad(Jx)$, which implies $\Tr(\ad(x) J)=-\Tr(\ad(Jx))$. This identity together with the condition $\Tr(\ad(x) J)=\Tr(\ad(Jx))$, which comes from $\psi\equiv 0$, and the fact that $J$ is an isomorphism imply the result.
	\end{proof}
	
	In dimension 6, there is only one unimodular non-nilpotent Lie algebra admitting an abelian complex structure (see \cite{ABD}). It is the Lie algebra $\mathfrak{s}$ determined by a basis $\{e_i\}_{i=1}^6$ and Lie brackets
	\begin{gather*} 
		[e_1,e_6]=-e_1, \quad [e_2,e_5]=e_1, \quad [e_1,e_5]=-e_2,\quad [e_2,e_6]=-e_2,\\
		[e_3,e_6]=e_3,\quad [e_4,e_5]=-e_3,\quad [e_3,e_5]=e_4,\quad [e_4,e_6]=e_4.
	\end{gather*}
	This Lie algebra appears as $\g_8$ in \cite{FOU} and as $\s_{(-1,0)}$ in \cite{ABD}; it is the real Lie algebra underlying the complex parallelizable Nakamura manifold \cite{Nak} (see also \cite{Ang}). It has an infinite number of non-equivalent\footnote{Two complex structures $J_1,J_2$ on a Lie algebra $\g$ are said to be equivalent if there exists a Lie algebra isomorphism $\varphi:\g\to\g$ such that $\varphi\circ J_1=J_2\circ \varphi$.} complex structures admitting a nonzero holomorphic $(3,0)$-form but only one of them is abelian (see \cite[Proposition 3.7]{FOU}), namely: $Je_1=e_2$, $Je_3=e_4$ and $Je_5=e_6$. It was proven in \cite{Yamada} that its corresponding simply connected Lie group  $S$ admits a lattice. We show next that this example can be generalized to any dimension of the form $4n+2$.
	
	\begin{example}\label{ex: nakamura-general}
		For $n\geq1$, let $\s_n=\R^2\ltimes \R^{4n}$ be the $(4n+2)$-dimensional unimodular Lie algebra with basis $\{f_1, f_2, e_1, e_2, \ldots, e_{4n}\}$ and Lie brackets given by\footnote{Throughout the article we use $A\oplus B$ to denote the block-diagonal matrix $\left[\begin{smallmatrix}A&\\&B\end{smallmatrix}\right]$. This naturally generalizes to the sum of $n$ square matrices.} \[  A:=\ad f_1|_{\R^{4n}}=\left(\matriz{0&-1\\1&0}\oplus \matriz{0&1\\-1&0}\right)^{\oplus n},\quad B:=\ad f_2|_{\R^{4n}}=\operatorname{diag}(1,1,-1,-1)^{\oplus n}.\] 
		Note that $\s_1$ coincides with the Lie algebra $\s$ above.
		
		It is easy to verify that the almost complex structure $J$ given by $Jf_1=f_2$ and $Je_{2j-1}=e_{2j}$ for all $1\leq j \leq 2n$ is abelian. It follows from Corollary \ref{corollary: abelianFCT} that $\s_n$ admits a nonzero closed $(n,0)$-form. We show next that the corresponding simply connected Lie group $S_n$ admits lattices. For $m\in\N$, $m\geq 3$, let $t_m=\log(\frac{m+\sqrt{m^2-4}}{2})$. Then
		\[ \exp(\pi A)=-\I_{4n}, \quad \exp(t_m B)=\operatorname{diag}(\e^{t_m}, \e^{t_m}, \e^{-t_m}, \e^{-t_m})^{\oplus n}.\] 
		Using that $\e^{t_m}+\e^{-t_m}=m$, it is easily seen that there exists $P\in \operatorname{GL}(4n,\R)$ such that $P^{-1} \exp(t_m B) P=\left[\begin{smallmatrix}0&-1\\1&m\end{smallmatrix}\right]^{\oplus 2n}$, and it is clear that $P^{-1}(-\I_{4n})P=-\I_{4n}$, so the matrices $\exp(\pi A)$ and $\exp(t_m B)$ are simultaneously conjugate to integer unimodular matrices. According to Theorem \ref{theorem: yamada}, since any basis of $\R^{4n}$ is rational, the subgroup $\Gamma_m^n:=(\pi \Z \oplus t_m \Z ) \ltimes P\Z^4$ is a lattice of $S_n$. The complex solvmanifold $(\Gamma_m^n\backslash S_n, J)$ has trivial canonical bundle for any $m$.
	\end{example}
	
	\begin{remark}
		The Lie algebra $\mathfrak{s}_n$ also carries a bi-invariant complex structure $\tilde{J}$ given by \[\tilde{J}f_1=-f_2,\qquad \tilde{J}e_{2j-1}=e_{2j}, \quad 1\leq j\leq {2n},\] so that the solvmanifold $(\Gamma\backslash S_n, \tilde{J})$ is complex parallelizable for any lattice $\Gamma\subset G$, generalizing in this way the complex parallelizable Nakamura manifold. 
	\end{remark} 
	
	\smallskip
	
	\subsection{Examples in almost nilpotent solvmanifolds}
	
	In this subsection we exhibit some examples of complex almost nilpotent solvmanifolds admitting invariant trivializing sections of the canonical bundle. In particular we show how to use Theorem \ref{theorem: yamada} in the case of a non-abelian nilradical.
	
	The examples we provide are considered in \cite{FP2}, where they characterize two types of Hermitian structures on almost nilpotent Lie algebras whose nilradical has one-dimensional commutator ideal, that is, $\n=\h_{2\l+1}\oplus \R^h$ with $l,h\in \N$, where $\h_{2\l+1}$ is a Heisenberg Lie algebra. Recall that $\h_{2\l+1}=\text{span}\{e_1,\ldots,e_{2\l+1}\}$ with $e_1$ central and $[e_j,e_{\l+j}]=e_1$ for $2\leq j\leq \l+1$.
	
	\begin{example}\label{exs: AN 1}
		Let $\g_n=\R e_{4n+2}\ltimes_B \h_{4n+1}$ where $B=\left[\begin{array}{c|cc} 
			0& & \\
			\hline 
			& \I_{2n} & \\
			&   & -\I_{2n} \end{array}
		\right]$ in the ordered basis $\{e_1,\ldots,e_{4n+1}\}$. 
		According to \cite[Proposition 2.4]{FP2}, the complex structure $J$ on $\g_n$ defined by $Je_1=e_{4n+2}$ and $Je_{2k}=e_{2k+1}$, $1\leq k\leq 2n$, is integrable. Moreover, it is easily verified that the associated Koszul form vanishes so that $(\g_n,J)$ admits a nonzero closed $(2n+1,0)$-form. For any $m\in \N$, $m\geq 3$, the associated simply connected Lie group $G_n$ admits a lattice $\Gamma_m^n$. Indeed, for $t_m=\log(\frac{m+\sqrt{m^2-4}}{2})$, let 
		\[ P_m=\left[\begin{array}{c|cc}
			1&0&0\\
			\hline 
			0& \I_{2n} & \alpha_m \I_{2n} \\ 0 & \frac{1}{\alpha_m^{-1}-\alpha_m} \I_{2n} & \frac{\alpha_m^{-1}}{\alpha_m^{-1}-\alpha_m} \I_{2n}  
		\end{array} \right], \quad \text{where}\quad \alpha_m=\exp(t_m).\] 
		Then $P_m^{-1} \exp(t_m B) P_m=\left[\begin{array}{c|cc} 1 & 0 &0 \\ \hline 
			0 & 0_{2n} &-\I_{2n}\\ 0& \I_{2n} & m \I_{2n} \end{array}\right]$. Thus, if we set $f_j= P_m e_j$, $1\leq j\leq 4n+1$, then we have that $[f_j, f_k]=[e_j,e_k]$, $1\leq j,k\leq 4n+1$, hence $\{f_j\}_{j=1}^{4n+1}$ is a rational basis of $\h_{4n+1}$ in which the matrix of $\exp(t_m B)$ is an integer unimodular matrix. It follows from Theorem \ref{theorem: yamada} that $\Gamma_m^n=t_m \Z\oplus \exp^{H_{4n+1}} (\text{span}_\Z \{f_1,\ldots,f_{4n+1}\})$ is a lattice of $G_m$. All the complex solvmanifolds $(\Gamma_m^n \backslash G_m,J)$ have trivial canonical bundle. 
	\end{example}
	
	\smallskip
	
	\begin{example}
		For $a_1,\ldots, a_n\in \R$ let us define $\g:=\g(a_1,\ldots,a_n)= \R e_{2n+2}\ltimes_B \h_{2n+1}$ where \[B=\left[\begin{array}{c|cc} 
			0 & & \\ \hline
			& 0 &-X \\
			& X & 0\end{array}
		\right], \quad X=\operatorname{diag}(a_1,\ldots,a_n),\] in the ordered basis $\{e_1,\ldots,e_{4n+1}\}$. The complex structure $J$ defined by $J e_1=e_{2n+2}$ and $Je_{2k}=e_{2k+1}$, $1\leq k\leq n$ is integrable, again due to \cite[Proposition 2.4]{FP2}. It is easy to verify that the Koszul form vanishes if and only if $\sum_{j=1}^n a_j=0$, and in this case $(\g,J)$ has a nonzero closed $(n+1,0)$-form. The simply connected Lie group $G:=G(a_1,\ldots,a_n)$ admits a lattice for some values of the parameters $a_1,\ldots,a_n$. Indeed, for any $n\in \N$ one can choose
		$a_1,\ldots,a_n\in\{2\pi,\pi,\frac{\pi}{2}\}+2\pi \Z$ with $a_1+\cdots+a_n=0$, and then $\{e_j\}_{j=1}^{2n+1}$ is a rational basis of $\h_{2n+1}$ in which $\exp B$ is a unimodular integer matrix, so by Theorem \ref{theorem: yamada} the Lie group $G$ admits a lattice $\Gamma:=\Gamma(a_1,\ldots,a_n)$. Then $(\Gamma\backslash G,J)$ has trivial canonical bundle.
	\end{example}

	\medskip
	
	\subsection{Holomorphically torsion canonical bundle}
	
	Here we consider the case when some power of the canonical bundle of a compact complex quotient $(\Gamma\backslash G,J)$ is trivialized by an invariant holomorphic section. We obtain that actually the canonical bundle is itself trivial by an invariant holomorphic section.
	
	\begin{proposition}
		Let $(G,J)$ be a $2n$-dimensional Lie group equipped with a left invariant complex structure. If $K_{(G,J)}^{\otimes k}$ admits a nonzero invariant holomorphic section for some $k\in \N$ then $K_{(G,J)}$ admits a nonzero invariant holomorphic section. That is, $(G,J)$ has trivial canonical bundle. The same happens for any quotient $\Gamma\backslash G$ where $\Gamma$ is a uniform lattice of $G$.
	\end{proposition}
	
	\begin{proof}
		We can work at the Lie algebra level since we are dealing with invariant objects. Let $\sigma$ be a generator of $\alt^{n,0} \g^*$, where $\g=\operatorname{Lie}(G)$. Then, $\sigma^{\otimes k}:=\sigma \otimes \cdots \otimes \sigma$ ($k$ times) is a generator of $(\alt^{n,0} \g^*)^{\otimes k}$, which we may assume holomorphic since this space is 1-dimensional. Recall from Remark \ref{remark: dsigma} that $d\sigma=\beta \wedge \sigma$ for some $(0,1)$-form $\beta$, which in terms of the extended Dolbeault operator $\bar{\partial}$ from \S\ref{S: prelim}  can be expressed as $\bar{\partial}\sigma=\beta \otimes \sigma$. 
		Next we compute 
		\[ 0=\bar{\partial} \sigma^{\otimes k}=\sum_{j=1}^k \sigma \otimes \cdots \otimes \underbrace{\bar{\partial}\sigma}_{j\text{-th place}} \otimes \cdots \otimes \sigma=\sum_{j=1}^k \beta \otimes \sigma^{\otimes k}=k \beta \otimes \sigma^{\otimes k}.\] Therefore, $\beta=0$ and this implies $\bar{\partial}\sigma=0$. Hence, $\sigma$ is holomorphic and the proof follows.
	\end{proof}

	\medskip
	
	\section{Triviality of the canonical bundle of solvable Lie groups with left invariant complex structures} \label{S: groups}
	
	The main goal in this section is to show that any simply connected solvable Lie group equipped with a left invariant complex structure has trivial canonical bundle. In general the trivializing holomorphic section will not be left invariant. 
	
	Any nowhere vanishing section of the canonical bundle of a $2n$-dimensional Lie group $G$ equipped with a left invariant complex structure $J$ can be written as $\tau=f \sigma$, where $\sigma$ is a nonzero left invariant $(n,0)$-form and $f:G\to \C^\times=\C\setminus\{0\}$ is a smooth function. If $\tau$ is closed we cannot expect uniqueness of the function $f$ in general (in the non-compact setting) as the following result shows.
	
	\begin{lemma}\label{lemma: holomorfas}
		Let $G$ be a $2n$-dimensional Lie group equipped with a left invariant complex structure $J$, and let $\sigma$ denote a nonzero left invariant $(n,0)$-form on $G$. Assume that $\tau_1:=f_1 \sigma$ is closed, for some smooth function $f_1:G\to \C^\times$. If $f_2:G\to \C^\times$ is another  smooth function on $G$ then $\tau_2:=f_2 \sigma$ is closed if and only if $\frac{f_2}{f_1}$ is a holomorphic function on $G$. 
		
		In particular, if $G$ is compact then $f_2=cf_1$ for some $c\in \C^\times$.
	\end{lemma}
	
	\begin{proof}
		Assume first that $H:=\frac{f_2}{f_1}$ is holomorphic. Then:
		\begin{align*}
			\overline{\partial} \tau_2 & = \overline{\partial} ((Hf_1)\sigma) = \overline{\partial}(Hf_1)\wedge \sigma + (Hf_1)\, \overline{\partial} \sigma = H (\overline{\partial}f_1) \wedge \sigma +(Hf_1)\,\overline{\partial} \sigma\\
			&= H(\overline{\partial}f_1 \wedge \sigma +f_1\,\overline{\partial} \sigma)= H\, \overline{\partial}(f_1\sigma) =0.
		\end{align*}
		Therefore $\tau_2$ is holomorphic and hence closed.
		
		Conversely, assume now that $\tau_2$ is closed. From $d\tau_1=0$ and $d\tau_2=0$ we obtain 
		\[ df_1\wedge \sigma + f_1\, d\sigma=0, \qquad df_2\wedge \sigma + f_2\, d\sigma=0.\]
		From these equations we obtain readily that
		\begin{equation}\label{eq:F2F1} d\left( \frac{f_2}{f_1}\right)\wedge \sigma=0.
		\end{equation}
		Let us consider a basis $\{\gamma_1,\ldots, \gamma_n\}$ of left invariant $(1,0)$-forms on $G$, hence we may assume $\sigma=\gamma_1\wedge\cdots\wedge \gamma_n$. If we write
		\[ d\left( \frac{f_2}{f_1}\right)=\sum_{j=1}^n (a_j\gamma_j+b_j\overline{\gamma}_j) \]
		for some $a_j,b_j\in\C$, then \eqref{eq:F2F1} becomes
		\[
		0 = \sum_{j=1}^n (a_j\gamma_j+b_j\overline{\gamma}_j)\wedge \sigma = \sum_{j=1}^n b_j\overline{\gamma}_j\wedge \sigma
		\]
		which implies $b_j=0$ for $j=1,\ldots,n$ since $\{\overline{\gamma}_j\wedge \sigma\}_{j=1}$ is a linearly independent set. This means that $d\left(\frac{f_2}{f_1}\right)$ is a $(1,0)$-form, and this is equivalent to $\overline{\partial}\left(\frac{f_2}{f_1}\right)=0$, that is, $\frac{f_2}{f_1}$ is holomorphic.
	\end{proof}
	
	\begin{corollary}\label{corollary: cpt uniqueness}
		With notation as in Lemma \ref{lemma: holomorfas}, if $\sigma$ is a nonzero invariant $(n,0$)-form on $\Gamma\backslash G$ and $f:\Gamma\backslash G\to \C^\times$ is a smooth function such that $\tau:=f\sigma$ is closed then $f$ is unique up to a nonzero constant. In particular, nowhere vanishing closed $(n,0)$-forms $\tau$ on $(\Gamma\backslash G,J)$ are either all invariant or all non-invariant.
	\end{corollary}

	
	\medskip
	
	Now we proceed to prove the main theorem of the section. We begin with a series of preliminary results. Recall that in a solvable Lie algebra $\g$ its nilradical $\n(\g)$ is given by $\n(\g)=\{x\in \g\mid \ad x \, \text{ is nilpotent}\}$.
	
	\begin{lemma}\label{lemma:cap}
		Let $\g$ be a solvable Lie algebra equipped with a complex structure $J$, and denote $\n(\g)$ its nilradical. If $\h=\Ker \psi$, where $\psi$ is the Koszul 1-form on $(\g,J)$ then \[\n(\g)\cap J\n(\g)\subseteq \h\cap J\h.\]
	\end{lemma}
	
	\begin{proof}
		Let $x\in \n(\g)\cap J\n(\g)$. Since $Jx\in \n(\g)$ then $\ad (Jx)$ is nilpotent and thus $\Tr \ad(Jx)=0$. As a consequence we only need to prove that $\Tr(J\ad x)=0$.
		It follows that 
		\[ x-iJx\in \n(\g)\oplus i\n(\g)=\n(\g)_\C=\n(\g_\C), \]
		so that $\ad(x-iJx)$ is a nilpotent endomorphism of $\g_\C$. We can write 
		\[ \ad(x-iJx)=\left[\begin{array}{c|c} 
			A_x & * \\
			\hline 
			0 & B_x
		\end{array}\right],\]
		in a certain basis of $\g_\C$ adapted to the decomposition $\g_\C=\g^{1,0}\oplus \g^{0,1}$.
		Since this operator is nilpotent, we have that both matrices $A_x$ and $B_x$ are nilpotent, so that $\Tr A_x=\Tr B_x=0$. Now we compute
		\[ J^\C \ad(x-iJx)=  \left[\begin{array}{c|c} 
			i A_x & * \\
			\hline 
			0 & -i B_x
		\end{array}\right]. \]
		Therefore, 
		\[ 0 = \Tr(J^\C \ad(x-iJx))=\Tr(J\ad(x))-i\Tr(J\ad(Jx)), \]
		so that 
		\[ \Tr(J\ad(x))=\Tr(J\ad(Jx))=0,\]
		that is, $x\in \h\cap J\h$.
	\end{proof}
	
	\medskip
	
	The following technical lemma will provide a particular basis of $(1,0)$-forms which will be useful in the proof of the main result of this section.
	
	\begin{lemma}\label{lemma:nice-basis}
		Let $\g$ be a $2n$-dimensional solvable Lie algebra equipped with a complex structure $J$. Then there exists a basis $\{\gamma_1,\ldots,\gamma_n\}$ of $(1,0)$-forms, with $\gamma_k=u^k+iv^k$, and an index $1\leq s\leq n$ such that:
		\begin{enumerate} 
			\item[$\ri$] $u^j$ is closed for $1\leq j\leq s$, and 
			\item[$\rii$] $u_j,v_j\in[\g,\g]\cap J[\g,\g]$  for $j>s$,
		\end{enumerate}
		where $\{u_1,v_1,\ldots,u_n,v_n\}$ denotes the dual basis of $\{u^1,v^1,\ldots,u^n,v^n\}$.
	\end{lemma}
	
	\begin{proof}
		Consider the commutator ideal $\g'=[\g,\g]$ and let $\u$ be a complementary subspace to $\g' \cap J\g'$ in $\g'$, that is \[\g'=(\g'\cap J\g') \oplus \u. \] Moreover, we have that $\g' \cap J\u=\{0\}$. Indeed, if $v\in \g'\cap J\u$, this implies that $v\in \g'\cap J\g'$ since $\u\subset \g'$. Hence, $Jv\in \u \cap (\g'\cap J\g')$, which implies $Jv=0$ and thus $v=0$.
		
		Therefore 
		we can decompose $\g$ as 
		\[\g=(\g'\cap J\g') \oplus \u \oplus J\u \oplus \v,\] where $\v$ is a complementary subspace to $\g' \oplus J\u$ in $\g$, which can be chosen $J$-invariant. As $\g$ is solvable, $\g'$ is a proper subspace of  $\g$, so $J\u\oplus \v\neq \{0\}$. The fact that the subspaces $\v$, $\u\oplus J\u$ and $\g'\cap J\g'$ are $J$-invariant allows us to take bases $\{x_1, \ldots,x_r, \tilde{x}_1,\ldots,\tilde{x}_r\}$ of $\v$, $\{y_1,\ldots,y_m, \tilde{y}_1,\ldots,\tilde{y}_m\}$ of $\u\oplus J\u$ (with $y_k\in J\u, \tilde{y}_k\in \u$) and $\{z_1,\ldots,z_\l, \tilde{z}_1,\ldots,\tilde{z}_\l\}$ of $\g'\cap J\g'$ such that $Jx_k=\tilde{x}_k$, $Jy_k=\tilde{y}_k$, $Jz_k=\tilde{z}_k$ and $r+m+\l=n$. Then, we can take the ordered basis of $(1,0)$-forms \[\{\gamma_1,\ldots,\gamma_s,\ldots,\gamma_n\}= \{x^1+i\tilde{x}^1,\ldots, x^r+i\tilde{x}^r, y^1+i\tilde{y}^1,\ldots, y^m+i\tilde{y}^m, z^1+i\tilde{z}^1,\ldots,z^\l+i\tilde{z}^\l\},\] with $s:=r+m<n$ and $\{x^1,\tilde{x}^1,\ldots,y^1,\tilde{y}^1,\ldots,z^1,\tilde{z}^1,\ldots\}$ is the basis of $\g^*$ dual to the basis $\{x_1,\tilde{x}_1,\ldots,y_1,\tilde{y}_1,\ldots,z_1,\tilde{z}_1,\ldots\}$. Let us rename $u^j=\operatorname{Re}\gamma_j$ and $v^j=\operatorname{Im}\gamma_j$. For $1\leq j\leq s$ we have that $u^j$ belongs to the annihilator of $\g'$ so $u^j$ is closed, and for $j>s$ we have that $u_j,v_j\in \g'\cap J\g'$. 
	\end{proof}
	
	\smallskip
	
	\begin{remark}
		It follows from the proof of Lemma \ref{lemma:nice-basis} that $s=n$ if and only if $\g=\g'\oplus J\g'$. In this case the complex structure $J$ is abelian. Indeed, the more general condition $\g'\cap J\g'=\{0\}$ implies that $J$ is abelian, which can be easily verified from $N_J=0$.
	\end{remark}
	
	\medskip
	
	\begin{theorem}\label{theorem:non-invariant}
		Any $2n$-dimensional simply connected solvable Lie group $G$ equipped with a left invariant complex structure $J$ admits a nonzero closed $(n,0)$-form $\tau$. In particular, the canonical bundle of $(G,J)$ is trivial.
	\end{theorem}
	
	\begin{proof}
		Let $\g$ be the Lie algebra of $G$ and take the basis $\{\gamma_1,\ldots,\gamma_n\}$ with $\gamma_k=u^k+iv^k$, as in Lemma \ref{lemma:nice-basis}.  
		Consider now the $(n,0)$-form $\sigma$ given by $\sigma=\gamma_1 \wedge \cdots \wedge \gamma_n$. If $d\sigma=0$ we may simply choose $\tau=\sigma$. On the other hand, if $d\sigma\neq 0$ it follows from Remark \ref{remark: dsigma} that
		\begin{equation}\label{eq:Cj}
			d\sigma=\frac14 \sum_{j=1}^n (-\psi(v_j)+i\psi(u_j)) \, \bar{\gamma_j} \wedge \sigma,
		\end{equation}
		Let us call $C_j=-\psi(v_j)+i\psi(u_j)$. We show first that $C_j=0$ when $j>s$. Indeed, in this case  $u_j,v_j\in \g'\cap J\g'$. As $\g'\cap J\g'\subset \n(\g)\cap J\n(\g)$ since $\g$ is solvable, it follows from Lemma \ref{lemma:cap} that $\psi(u_j)=\psi(v_j)=0$ so $C_j=0$ for $j>s$. Therefore we can write
		\begin{align*} 
			d\sigma=\frac14 \sum_{j=1}^s C_j \bar{\gamma_j} \wedge \sigma
			=\frac14 \sum_{j=1}^s (C_j \bar{\gamma_j} \wedge \sigma+ C_j \underbrace{\gamma_j \wedge \sigma}_{=0})
			=\frac12 \sum_{j=1}^s C_j u^j \wedge \sigma.
		\end{align*}
		Hence, the form $\alpha=\frac12 \sum_{j=1}^s C_j u^j$ satisfies $d\sigma=\alpha\wedge \sigma$ and $d\alpha=0$ by the choice of the basis $\{\gamma_1,\ldots,\gamma_n\}$. Since $G$ is simply connected the left invariant 1-form $\alpha$ on $G$ is exact so that there exists a smooth function $f:G\to \C$ satisfying $\alpha=df$. Finally, we consider the $(n,0)$-form $\tau:=\e^{-f} \sigma$ and we compute \[ d\tau=\e^{-f} (-\alpha \wedge \sigma+d\sigma)=0,\] which says that $\tau$ is a nowhere vanishing closed $(n,0)$-form on $G$. 
	\end{proof}
	
	\smallskip 
	
	\begin{remark}\label{remark:F hom}
		The closed $(n,0)$-form $\tau$ from Theorem \ref{theorem:non-invariant} can be written as $\tau=F \sigma$, where $F:G\to\C^\times$ is a Lie group homomorphism. Indeed, with the notation in the proof of Theorem \ref{theorem:non-invariant}, replacing $f$ by $f-f(1_G)$ we still have that $\alpha=df$ is left invariant and $f(1_G)=0$, where $1_G$ denotes the identity element. This implies that $f:G\to\C$ is an additive homomorphism. Hence, $F:=\e^{-f}:G\to\C^\times$ is a multiplicative homomorphism.
	\end{remark}

	
	
	
	
	\smallskip
	
	\begin{remark}
		There are Lie groups with a left invariant complex structure which do not have trivial canonical bundle. For instance, the Hopf manifold $\mathbb{S}^1\times \mathbb{S}^3 \cong \mathbb{S}^1\times \operatorname{SU}(2)$ carries a left invariant complex structure, and this compact complex surface has non-trivial canonical bundle. 
	\end{remark}
	
	\smallskip
	
	\begin{remark}
		It was conjectured by Hasegawa in \cite{Has2} that all simply connected unimodular solvable Lie groups with left invariant complex structure are
		Stein manifolds (that is, they are biholomorphic to a closed complex submanifold of some $\C^{N}$). If this conjecture were true, then the canonical bundle of any of these pairs $(G,J)$ would be holomorphically trivial according to the Oka-Grauert principle (\cite{Grauert}), since the canonical bundle of any such Lie group is always smoothly trivial via a left invariant section. Thus, Theorem \ref{theorem:non-invariant} provides evidence in the direction of this conjecture.
	\end{remark}
	
	\smallskip
	
	\section{An algebraic obstruction for the triviality of the canonical bundle} \label{S: solvmflds}
	
	In this section we will consider compact complex manifolds obtained as quotients of a Lie group by a uniform lattice. We provide an algebraic obstruction for the canonical bundle to be holomorphically trivial (or more generally, holomorphically torsion), in terms of the Koszul 1-form $\psi$. Namely, $\psi$ has to vanish on the commutator ideal of the associated Lie algebra. We will do this by exploiting the relation of $\psi$ with the Chern-Ricci form of any invariant Hermitian metric on the quotient.
	
	Let us recall the definition of the Chern-Ricci form. Let $(M,J,g)$ be a $2n$-dimensional Hermitian manifold, and let $\omega=g(J\cdot,\cdot)$ be the fundamental 2-form associated to $(M,J,g)$. The \textit{Chern connection} is the unique connection $\nabla^C$ on $M$ which is Hermitian (i.e. $\nabla^C J=0$, $\nabla^C g=0$) and the $(1,1)$-component $T^{1,1}$ of its torsion tensor vanishes. In terms of the Levi-Civita connection $\nabla$ of $g$, the Chern connection is expressed as
	\[ g(\nabla^C_X Y,Z)=g(\nabla_X Y,Z)-\tfrac12 d\omega(JX,Y,Z),\quad X,Y,Z\in\X(M).\]
	The \textit{Chern-Ricci form} $\rho=\rho(J,g)$ is defined by 
	\[ \rho(X,Y)=-\tfrac12 \Tr( J \circ R^C(X,Y))=\sum_{i=1}^n g(R^C(X,Y)e_i,J e_i),\] where $R^C(X,Y)=[\nabla^C_X, \nabla^C_Y]-\nabla^C_{[X,Y]}$ is the curvature tensor associated to $\nabla^C$ and $\{e_i,Je_i\}_{i=1}^n$ is a local orthonormal frame for $g$. It is well known that $\rho$ is a closed real $(1,1)$-form on $M$.
	
	\smallskip
	
	Consider now any left invariant almost Hermitian structure $(J,g)$ on a simply connected Lie group  $G$, not necessary solvable, with Lie algebra $\g$. In \cite{V} it is proved that 
	\begin{equation}\label{eq: Chern-Ricci form}
		\rho(x,y)=\tfrac{1}{2}(\Tr(J \ad[x,y])-\Tr \ad(J [x,y])), \quad x,y\in\g.
	\end{equation}
	Remarkably, this Chern-Ricci form does not depend on the Hermitian metric. We observe from \eqref{eq: Chern-Ricci form} that $2\rho=-d\psi$, thus implying that if $\psi$ vanishes then $\rho=0$, for any Hermitian metric $g$. 
	As a consequence we have:

	\begin{proposition}\label{proposition: chern-ricci}
		If there exists a nonzero left invariant holomorphic $(n,0)$-form on $G$ then for any left invariant Hermitian metric $g$ on $G$, the induced Hermitian structure $(J,g)$ on $\Gamma\backslash G$ has vanishing Chern-Ricci form. In particular, the restricted Chern holonomy of $(J,g)$ on $\Gamma\backslash G$ is contained in $\operatorname{SU}(n)$. 
	\end{proposition}
	
	\begin{proof}
		We only have to justify the last statement, and this follows from \cite[Proposition 1.1]{To}.
	\end{proof}
	
	\smallskip
	
	Thus, if $\rho\neq 0$ (that is, $\psi([\g,\g])\neq 0$) then there is no invariant trivializing section of the canonical bundle of $(\Gamma\backslash G,J)$. We show next that the condition $\rho \neq 0$ is also sufficient to prove that the canonical bundle is not holomorphically trivial. In fact, we will prove a stronger result, namely that if $\rho\neq 0$ then the canonical bundle of $(\Gamma\backslash G,J)$ is not holomorphically torsion. 
	
	Following ideas from \cite[Proposition 5.1]{GP}, our result will be proved using Belgun's symmetrization, which we state below.
	
	\begin{lemma}\label{lemma: sym}$($\cite[Theorem 7]{Belgun}, \cite[Theorem 2.1]{FG}$)$
		Let $M=\Gamma\backslash G$ be a compact quotient of a simply connected Lie group by a uniform lattice $\Gamma$ with an invariant complex structure $J$. Let $\nu$ denote the bi-invariant volume form on $G$ given in \cite[Lemma 6.2]{Mil} and such that $\int_M \nu=1$. Identifying left invariant forms on $M$ with linear forms over $\g^*$ via left translations, consider the Belgun symmetrization map defined by:
		\[ \mu :\Omega^*(M)\to \alt^* \g^*, \quad \mu(\alpha)(X_1,\ldots,X_k)=\int_M \alpha_m(X_1|_m,\ldots, X_k|_m)\nu_m, \]
		for $X_1,\ldots,X_k\in\mathfrak{X}(M)$. Then:
		\begin{enumerate}
			\item[$\ri$] $\mu(f)\in\R$ for any $f\in C^\infty(\Gamma\backslash G,\R)$,  
			\item[$\rii$] $\mu(\alpha)=\alpha$ if $\alpha\in \alt^* \g^*$;
			\item[$\riii$] $\mu(J\alpha)=J\mu(\alpha)$, where $J\alpha(\cdot,\ldots,\cdot)=\alpha(J^{-1}\cdot,\ldots,J^{-1}\cdot)$;
			\item[$\riv$] $\mu(d\alpha)=d(\mu(\alpha))$.
		\end{enumerate}
		Extending this map $\C$-linearly to $\C$-valued differential forms on $M$, we also have:
		\begin{enumerate}
			\item[$\rm{(v)}$] $\mu(\partial\alpha)=\partial(\mu(\alpha))$ and $\mu(\overline{\partial}\alpha)=\overline{\partial}(\mu(\alpha))$.
		\end{enumerate}
	\end{lemma}

	\smallskip
	
	\begin{theorem}\label{theorem: obstruction}
		If the canonical bundle of $(\Gamma\backslash G,J)$ is trivial (or, more generally, holomorphically torsion) then $\Tr (J \ad([x,y]))=0$ for all $x,y\in \g=\operatorname{Lie}(G)$, that is $\psi([\g,\g])= 0$.
	\end{theorem}
	
	\begin{proof}
		According to \cite[Proposition 1.1]{To} if $M$ is a compact complex manifold and $K_M$ is holomorphically torsion  then given any Hermitian metric $g$ on $M$, the associated Chern-Ricci form $\rho$ satisfies $\rho=i\partial\bar{\partial}F$, for some $F\in C^{\infty}(M,\R)$.   
		
		Now, assume that the canonical bundle of $(\Gamma\backslash G,J)$ is holomorphically torsion and consider a Hermitian metric $g$ on $\Gamma\backslash G$ induced by a left invariant one on $G$. Then its associated Chern-Ricci form $\rho$ satisfies $\rho=i\partial\bar{\partial}F$, for some $F\in C^\infty(\Gamma\backslash G,\R)$. We consider next the symmetrization $\mu(\rho)$ of $\rho$. It follows from Lemma \ref{lemma: sym}(v) that 
		\[ \mu(\rho)=i\partial\bar{\partial}\mu(F)=0,\] 
		since $\mu(F)$ is constant. As $\rho$ is left invariant we obtain $\rho=\mu(\rho)$ and therefore $\rho=0$. Since at the Lie algebra level  we have that $\rho(x,y)=\Tr(J\ad([x,y]))$ for $x,y\in\g$, the proof is complete.
	\end{proof}
	
	\smallskip
	
	\begin{remark}
		Proposition 1.1 in \cite{To} predicts the existence of a Hermitian metric with $\rho=0$ on any compact complex manifold with holomorphically torsion canonical bundle. It follows from the proof of Theorem \ref{theorem: obstruction} that in the case of a complex compact quotient $\Gamma\backslash G$ any \textit{invariant} Hermitian metric has $\rho=0$. Using again \cite[Proposition 1.1]{To}, we obtain that these invariant Hermitian metrics have restricted Chern holonomy contained in $\SU(n)$, where $2n$ is the real dimension of the manifold.
	\end{remark}
	
	\begin{remark}
		If $\rho=0$ then the canonical bundle of the complex solvmanifold is not necessarily trivial. Indeed, consider the Lie group $G$ from Example \ref{ex: motivation}  equipped with the left invariant complex structure $J$ given therein. It is easy to see that $\psi(e_0)=-2$ and $\psi(e_j)=0$ for $1\leq j\leq 3$, so that $\psi([\g,\g])=0$. However, $G$ admits a lattice $\Gamma'=\{(\pi k, m, n, \frac{p}{2}) \mid k,m,n,p\in \Z\}$ such that $(\Gamma'\backslash G, J)$ is a secondary Kodaira surface (see \cite{Hase}) and hence, has non-trivial canonical bundle. Note that $\tau \otimes \tau$, where $\tau=\e^{it} \sigma$, is a trivializing section of $K_{(\Gamma'\backslash G,J)}^{\otimes 2}$, and thus the canonical bundle is holomorphically torsion. 
		
		\medskip
		
		We exhibit next a 6-dimensional example of this phenomenon. 
		
		\begin{example}
			For $p\in \R$, let $\g_p=\R e_6\ltimes_{A_p} \R^5$, where the matrix $A_p$ is given in the basis $\{e_1,\ldots,e_5\}$ of $\R^5$ by: 
			\[A_p=\matriz{-p&-1&&&\\1&-p&&&\\&&p&2&\\&&-2&p&\\&&&&0}.\]
			Equip $\g_p$ with the complex structure $Je_1=e_2$, $Je_3=e_4$ and $Je_5=e_6$. Then an easy calculation shows that $\psi(e_j)=0$ for $1\leq j\leq 5$ and $\psi(e_6)=2$. It follows from Theorem \ref{theorem: inv} that this Lie algebra does not admit any nonzero closed $(3,0)$-form.
			
			For some values of $p\in \R$, the associated simply connected Lie group $G_p$ admits lattices, according to \cite{CM} (the Lie algebra $\g_p$ corresponds to the Lie group denoted by $G_{5.17}^{p,-p,2}\times \R$ there). Moreover, it was shown, using techniques by Console and Fino in \cite{CF}, that for certain values of $p$ some lattices $\Gamma$ in $G_p$ satisfy $b_3(\Gamma\backslash G_p)=0$ (see \cite[Table 7.1]{CM}); for instance, for any $m\in \N$ take $p=\frac{s_m}{\pi}$, where $s_m=\log\left(\frac{m+\sqrt{m^2+4}}{2}\right)$. Then, $\exp(\pi A_p)=\operatorname{diag}(-\e^{-s_m},-\e^{-s_m},\e^{s_m},\e^{s_m},1)$ is conjugate to the integer unimodular matrix $E_m=\matriz{0&1\\1&m}^{\oplus 2}\oplus (1)$. According to Theorem \ref{theorem: yamada} the Lie group $G_p$ admits a lattice $\Gamma_m=\pi \Z \ltimes P\Z^5$, where $P^{-1}\exp(\pi A_p)P=E_m$. Since $b_3(\Gamma_m\backslash G_p)=0$,  it follows from Proposition \ref{proposition:betti-n} that the canonical bundle of $(\Gamma_m\backslash G_p,J)$ is not holomorphically trivial. However, this complex solvmanifold has holomorphically torsion canonical bundle. Indeed, if $\sigma$ is a nonzero left invariant $(3,0)$-form on $G_p$ then $\tau \otimes \tau$, where $\tau=\e^{it}\sigma$, induces a trivializing section of $K_{\Gamma_m\backslash G_p}^{\otimes 2}$ since $\e^{2it}=(\e^{it})^2$ is $\pi$-periodic.
		\end{example}
	\end{remark}
	
	\smallskip
	
	\begin{example}
		There are examples of complex solvmanifolds whose canonical bundle is not holomorphically torsion (and in particular not holomorphically trivial). Such examples are given by Oeljeklaus-Toma manifolds, introduced in \cite{OT}.
		These complex manifolds were constructed from certain number fields, but later Kasuya showed in \cite{Ka} that they are complex solvmanifolds. 
	\end{example}
	
	\medskip
	
	As another illustration of the obstruction from Theorem \ref{theorem: obstruction}, we deal in the next result with the case of compact semisimple Lie groups. We recall Samelson's construction of a complex structure on a compact semisimple even-dimensional Lie algebra $\g$ \cite{Sam}.
	
	Let $\h$ be a maximal abelian subalgebra of $\g$. Then we have the root space decomposition of $\g_{\C}$ with respect to $\h_{\C}$
	\[ \g_{\C}= \h_{\C} \oplus \sum _{\alpha \in
		\Phi} \g_{\alpha},\] 
	where $\Phi$ is the finite subset of nonzero elements in 
	$(\h_{\C})^*$ called roots, and 
	\[\g_{\alpha}=\{ x \in \g_{\C} \mid  [h,x]=\alpha(h)x \quad \forall h \in  \h_{\C} \}\] are the
	one-dimensional root subspaces. Since $\h$ is even-dimensional, one can choose a skew-symmetric endomorphism $J_0$ of $\h$ with
	respect to the Killing form such that $J_0^2=-\I_\h$. Samelson defines a complex structure on $\g$ by considering a positive system
	$\Phi^+$ of roots, which is a set $\Phi^+ \subset \Phi$ satisfying
	\[ \Phi^+ \cap (-\Phi^+) = \emptyset, \qquad \Phi^+ \cup (-\Phi^+)
	=\Phi, \qquad \alpha , \beta \in \Phi ^+ , \; \alpha +\beta \in
	\Phi \Rightarrow \alpha +\beta \in \Phi ^+ .\] Setting \[ \m =
	\h^{1,0} \oplus \sum _{\alpha \in \Phi^+}  \g_{\alpha}, \] 
	where $\h^{1,0}$ is the eigenspace of $J_0^\C$ of eigenvalue $i$, it follows that $\m$ is a complex Lie subalgebra of
	$\g_{\C}$ which induces a complex structure $J$ on $\g$ such
	that $\g^{1,0}=\m$, that is, $\m$ is the eigenspace of $J^\C$ with eigenvalue $i$. This complex structure is skew-symmetric
	with respect to the Killing form on $\g$.
	
	Conversely, Pittie proved in \cite{Pit} that any left invariant complex structure on $G$ is obtained in this way.  
	
	\smallskip 
	
	In the next result we use Theorem \ref{theorem: obstruction} in order to show that the canonical bundle of a compact Lie group equipped with a left invariant complex structure is not holomorphically torsion.
	
	\begin{proposition}\label{proposition:compactos}
		The canonical bundle of a $2n$-dimensional compact semisimple Lie group equipped with a left invariant complex structure is not holomorphically torsion.
	\end{proposition}
	
	\begin{proof}
		We use the notation from the paragraphs above: $G$ is the compact Lie group, $\g$ its Lie algebra and $J:\g\to\g$ is the complex structure obtained by Samelson's construction. 
		
		Since $[\g,\g]=\g$,  according to Theorem \ref{theorem: obstruction} it is sufficient to show that $\Tr(J\ad x)\neq 0$ for some $x\in\g$, or equivalently, $\Tr(J^\C\ad x)\neq 0$ for some $x\in\g_\C$.
		
		Recall that $\g^{1,0}=\h^{1,0} \oplus \sum _{\alpha \in \Phi^+}  \g_{\alpha}$ and $\g^{0,1}=\h^{0,1} \oplus \sum _{\alpha \in \Phi^+}  \g_{-\alpha}$. Let $x_\alpha$ be a generator of $\g_\alpha$ for any $\alpha\in \Phi$. If $\{h_1,\ldots,h_r\}$ is a basis of $\h^{1,0}$, then
		$\mathcal{B}=\{h_1,\ldots,h_r\}\cup\{x_\alpha\mid \alpha\in\Phi^+\}$ is a basis of $\g^{1,0}$ and $\overline{\mathcal{B}}=\{\overline{h_1},\ldots,\overline{h_r}\}\cup\{x_{-\alpha}\mid \alpha\in\Phi^+\}$ is a basis of $\g^{0,1}$. 
		
		Consider now $h\in\h^{1,0}\subset \g^{1,0}$. Then, with respect to the basis $\mathcal{B}\cup  \overline{\mathcal{B}}$ of $\g$, we have:
		\[ \ad h=\left[\begin{array}{c|c} 
			A_h & * \\
			\hline 
			0 & B_h
		\end{array}\right] \qquad \text{and} \qquad J^\C \ad h=\left[\begin{array}{c|c} i A_h & * \\ \hline 
			0 &  -i B_h \end{array}\right].\]
		More precisely, since $\h^{1,0}$ is an abelian subalgebra and $[h,x_\alpha]=\alpha(h)x_\alpha$, 
		the matrices $A_h$ and $B_h$ are given by:
		\[ A_h=\left[\begin{array}{c|ccc}
			0_r & && \\
			\hline 
			& \alpha_1(h) & & \\
			&  &   \ddots & \\
			&  &     &    \alpha_s(h)
		\end{array} 
		\right] \quad \text{and} \quad  
		B_h=\left[\begin{array}{c|ccc}
			0_r & && \\
			\hline 
			& -\alpha_1(h) & & \\
			&  &   \ddots & \\
			&  &     &    -\alpha_s(h)
		\end{array} 
		\right],\]
		where $s=|\Phi^+|$. Hence, 
		\[ \Tr(J^\C \ad h)=2i\sum_{j=1}^s \alpha_j(h)=2i\sum_{\alpha\in\Phi^+} \alpha(h).   \]
		It is known that $\sum_{\alpha\in\Phi^+} \alpha\neq 0$. Indeed, there is $\Pi\subset \Phi^+$, whose elements are known as simple roots, such that $\Pi$ is a basis of $(\h_\C)^*$ and each $\alpha\in\Phi^+$ is a linear combination of the simple roots with non-negative integer coefficients. Therefore, if $\sum_{\alpha\in\Phi^+}\alpha =0$ then every $\alpha\in \Phi^+$ would be zero, which is impossible.
		
		As a consequence, we can choose $h\in\h^{1,0}$ such that $\Tr(J^\C\ad h)\neq 0$.
	\end{proof}
	
	\smallskip
	
	\begin{remark}
		Assume that $G$ is a non-compact semisimple Lie group equipped with a left invariant complex structure $J$. It follows from \cite{Bo} that $G$ has a uniform lattice $\Gamma$. Again, since $G$ is semisimple (so that $\g=\operatorname{Lie}(G)$ satisfies $[\g,\g]=\g$), it follows from Theorem \ref{theorem: inv} and Theorem \ref{theorem: obstruction} that the canonical bundle of the compact complex manifold $(\Gamma\backslash G,J)$ is either trivial via an invariant section (when $\psi= 0$) or it is not holomorphically torsion (when $\psi\neq 0$), where $\psi$ is the Koszul $1$-form on $(\g,J)$.
		
		Some recent results concerning non-compact semisimple Lie groups are the following:
		\begin{itemize}
			\item In \cite{GP} it was proved that any non-compact real simple Lie group $G$ of inner type and even dimension carries a left invariant complex structure $J$. Moreover, if $\Gamma$ is a lattice in $G$ then the canonical bundle of $(\Gamma\backslash G,J)$ is not holomorphically torsion.
			
			\item In \cite{OU} it was proved that a $6$-dimensional unimodular non-solvable Lie algebra admits a complex structure with a nonzero closed $(3,0)$-form if and only if it is isomorphic to $\mathfrak{so}(3, 1)$. It follows that $\Gamma\backslash\operatorname{SO}(3,1)$ carries a complex structure with trivial canonical bundle (via an invariant section) for any lattice $\Gamma$.
			
			As a generalization, we observe that if $\g$ is a  semisimple complex Lie algebra then its ``realification'' $\g_\R$ admits a bi-invariant complex structure $J$. If $G_\R$ denotes the simply connected Lie group associated to $\g_\R$ then the pair $(G_\R,J)$ admits a left invariant trivializing section of its canonical bundle, since it is complex parallelizable. Therefore, any compact quotient $(\Gamma\backslash G_\R,J)$ has trivial canonical bundle.
		\end{itemize}
	\end{remark}

	\subsection{More examples of complex solvmanifolds with trivial canonical bundle} We look for examples of complex solvmanifolds $(\Gamma\backslash G,J)$ with trivial canonical bundle when there are no invariant trivializing sections. Due to Theorem  \ref{theorem: obstruction} we need 
	$\psi\not\equiv 0$ but $\psi([\g,\g])=0$. We show next that in many cases 
	we obtain such a trivializing section.
	
	\begin{proposition}\label{proposition: tau explicit}
		Let $(G,J)$ be a $2n$-dimensional simply connected solvable unimodular Lie group with a left invariant complex structure. Let $\h$ denote the kernel of $\psi$ and assume that $\psi([\g,\g])\equiv 0$, so that $\g=\R e_0 \ltimes \h$ and consequently $G=\R\ltimes H$, where $H$ is the unique connected normal subgroup of $G$ such that $\operatorname{Lie}(H)=\h$. Then the $(n,0)$-form $\tau=\exp(-\frac{i}{2} \Tr (J \ad e_0) t) \sigma$ is closed, where $t$ is the coordinate of $\R$ and $\sigma$ is a left invariant $(n,0)$-form. 
	\end{proposition}
	
	\begin{proof}
		By using a Hermitian inner product on $\g$ we can choose $e_1\in \h\cap (\h\cap J\h)^\perp$, so that $\h=\R e_1\oplus (\h\cap J\h)$. Set next $e_0:=-Je_1\in \h^\perp$, hence $\g=\R e_0\ltimes \h$. Let $\{u_j,v_j\}_{j=1}^{n-1}$ be a basis of $\h\cap J\h$ such that $Ju_j=v_j$, $1\leq j\leq n-1$. 
		
		Define the $(n,0)$-form $\sigma$ on $\g$ by $\sigma=(e^0+ie^1) \wedge \gamma_1\wedge\cdots\wedge \gamma_{n-1}$, where $\gamma_j=u^j+iv^j$ and $\{e^0,e^1,u^i,v^i\}_{i=1}^{n-1}$ is the dual basis of $\{e_0,e_1,u_i,v_i\}_{i=1}^{n-1}$. In this basis,  Remark \ref{remark: dsigma} implies that 
		\begin{align*} d\sigma&= \tfrac{i}{4} \Tr (J \ad e_0)\, (e^0-ie^1) \wedge (e^0+ie^1) \wedge \gamma_1 \wedge \cdots \wedge \gamma_{n-1} \\
			&= -\tfrac12 \Tr(J \ad e_0)  e^{01} \wedge \gamma_1 \wedge \cdots \wedge \gamma_{n-1}.
		\end{align*}
		By the definition of the product on $G=\R\ltimes H$ it follows that if we consider the dual basis $\{e^0,e^1,u^i,v^i\}_{i=1}^{n-1}$ as left invariant 1-forms on $G$ (which is diffeomorphic to $\R^{2n}$), then $dt=e^0$. Let us consider the $(n,0)$-form given by $\tau=\e^{i\lambda t} \sigma$,  where $\lambda=-\frac12 \Tr(J \ad e_0)$. We compute
		\begin{align*}
			d\tau&=\e^{i\lambda t} \left(i\lambda\, e^0 \wedge \sigma +d\sigma\right)\\
			&=\e^{i\lambda t} \left(\tfrac12 \Tr(J\ad e_0)-\tfrac12 \Tr(J \ad e_0)  \right) e^{01}\wedge \gamma_1\cdots\wedge \gamma_{n-1}\\
			&=0.
		\end{align*}
		Therefore, $\tau$ is closed.
	\end{proof}
	
	
	\medskip 
	
	In the next example we apply Proposition \ref{proposition: tau explicit} in order to show the triviality of the canonical bundle associated to complex structures of splitting type (see \cite{AOUV} for a precise definition) on the 6-dimensional complex parallelizable Nakamura solvmanifold. 
	
	\begin{example}
		In \cite[Proposition 3.1]{AOUV} complex structures of splitting type on the 6-dimensional complex parallelizable Nakamura manifold are classified. There are three non-equivalent cases: 
		
		\begin{enumerate} 
			\item[\ri] $J:  d\omega^1=-\omega^{13}, \quad d\omega^2=\omega^{23}, \quad d\omega^3=0$,
			
			\medskip 
			
			\item[\rii] $J_A: \begin{cases} d\omega^1=A\omega^{13}-\omega^{1\bar{3}},\\
				d\omega^2=-A\omega^{23}+\omega^{2\bar{3}},\quad A\in\C, \; |A|\neq 1,\\
				d\omega^3=0,
			\end{cases}$
			
			\medskip 
			
			\item[\riii] $J_B: \begin{cases} d\omega^1=-\omega^{13}+B\omega^{1\bar{3}}, \\
				d\omega^2=-\bar{B}\omega^{23}+\omega^{2\bar{3}},\quad B\in\C,\; |B|<1, \\
				d\omega^3=0 
			\end{cases}$
		\end{enumerate}
		where $\{\omega^1,\omega^2,\omega^3\}$ is a basis of $(1,0)$-forms.
		
		According to \cite[Proposition 3.7]{FOU}, the underlying Lie algebra admits a nonzero holomorphic $(3,0)$-form only for complex structures of type (i) and (ii). Therefore any associated solvmanifold equipped with $J_B$ in (iii) admits no invariant trivializing sections of the canonical bundle. Nevertheless, we show next that there are lattices such that the associated complex solvmanifolds equipped with $J_B$ do have trivial canonical bundle.
		
		In a real basis of 1-forms $\{f^1,\ldots,f^6\}$ such that in the dual basis $J_B$ is given by $J_B f_{2i-1}=f_{2i}$, equations (iii)  can be written as
		\begin{align*}
			df^1 =  (r-1) f^{15} +s f^{16}- s f^{25} + (r+1) f^{26},\quad&
			df^2 =  s f^{15}-(r+1) f^{16}+(r-1)f^{25}+s f^{26},\\
			df^3 = (1-r) f^{35}-s f^{36}-s f^{45}+(r+1) f^{46},\quad& 
			df^4 = s f^{35}-(r+1)f^{36}-(r-1)f^{45}-sf^{46},\\
			df^5 =0,\quad& df^6=0
			,\end{align*} 
		where $B=r+is$. Therefore, the Lie brackets determined by $\{df^1,\ldots,df^6\}$ are
		\begin{align*} 
			[f_1,f_5]&=(1-r) f_1-s f_2,\quad [f_1,f_6]=-s f_1+(r+1) f_2,\\ [f_2,f_5]&=s f_1+(1-r) f_2, \quad [f_2,f_6]=-(r+1) f_1-s f_2, \\
			[f_3,f_5]&=(r-1) f_3-s f_4,\quad [f_3,f_6]=s f_3+(r+1) f_4,\\ [f_4,f_5]&=s f_3+(r-1) f_4,\quad [f_4,f_6]=-(r+1) f_3+s f_4.
		\end{align*} 
		
		Let us denote $\g$ the Lie algebra determined by these Lie brackets. On the other hand, recall the Lie algebra $\s$ from the paragraph before Example \ref{ex: nakamura-general}. It is straightforward to verify that $\varphi:(\g,J_B)\to (\s, \tilde{J_B})$ given by 
		\[\varphi=\matriz{0&0&0&1\\0&0&1&0\\0&-1&0&0\\1&0&0&0}\oplus\matriz{-s&r+1\\1-r&-s}\] is a biholomorphic isomorphism, where 
		\begin{gather*} \tilde{J_B} e_1=-e_2,\quad \tilde{J_B} e_2=e_1, \quad \tilde{J_B} e_3=e_4, \quad \tilde{J_B} e_4=-e_3,\\
			\tilde{J_B} e_5=\frac{-2s\, e_5+(r^2+s^2-2r+1)\, e_6}{r^2+s^2-1}, \quad  \tilde{J_B} e_6=\frac{-(r^2+s^2+2r+1)\,e_5+2s\, e_6}{r^2+s^2-1}.
		\end{gather*}
		
		The Koszul $1$-form on $(\s,\tilde{J_B})$ is given by $\psi=4e^5$. Since $\psi([\s,\s])=0$ we can apply Proposition \ref{proposition: tau explicit} and get a closed nowhere vanishing $(3,0)$-form in the Lie group $S=\R\ltimes H$ given by 
		\[\tau=\e^{-2it} (e^1-ie^2) \wedge (e^3+ie^4)\wedge \left(e^5-i\left(\frac{2s}{r^2+s^2-1} e^5+\frac{r^2+s^2+2r+1}{r^2+s^2-1}e^6\right)\right), \]
		where $t$ is the coordinate of $\R$.
		On the other hand, according to Example \ref{ex: nakamura-general}, the Lie group $S$ admits lattices given by $\Gamma_m=(\pi \Z\oplus t_m \Z) \ltimes P_m\Z^4$ for $m\in \N$, $m\geq 3$. The form $\tau$ is invariant by $\Gamma_m$ since $\exp(-2i(t+\pi k))=\exp(-2it)$ for all $k\in\Z$, so it induces a closed non-vanishing $(3,0)$-form on $(\Gamma_m\backslash S,\tilde{J_B})$, which therefore have trivial canonical bundle.
	\end{example} 
	
	We finish this section with an example of a solvmanifold with trivial canonical bundle such that its Lie algebra does not appear in \cite[Proposition 2.8]{FOU}.
	
	\begin{example}
		Let $\g=\R e_6 \ltimes_A \R^5$, where in the basis $\{e_1,\ldots,e_5\}$ the matrix $A$ is given by \[A:=\ad e_6|_{\R^5}=\left[\begin{smallmatrix}0&1&1&0&0\\-1&0&0&1&0\\0&0&0&1&0\\0&0&-1&0&0\\0&0&0&0&0 \end{smallmatrix}\right].\] 
		Equip $\g$ with the complex structure $Je_1=e_2, Je_3=e_4, Je_5=e_6$. By Theorem \ref{theorem: inv}, since $\psi=4e^6$, $(\g,J)$ does not admit a nonzero closed $(3,0)$-form so that $\g$ does not appear in \cite[Proposition 2.8]{FOU}. However, since $\psi([\g,\g])= 0$, using Proposition \ref{proposition: tau explicit} we obtain a closed non-vanishing $(3,0)$-form $\tau$ on the associated simply connected Lie group $G$, which is given by 
		$\tau=\e^{-2it} (e^1+ie^2)\wedge(e^3+ie^4)\wedge(e^5+ie^6)$, where $t$ is the coordinate of $\R$. Moreover, 
		\[\exp(\pi A)=\left[\begin{smallmatrix}-1&0&-\pi&0&0\\0&-1&0&-\pi&0\\0&0&-1&0&0\\0&0&0&-1&0\\0&0&0&0&1\end{smallmatrix}\right]\] is conjugate to its Jordan form $B:=\left[\begin{smallmatrix}
			-1&1 &0&0&0\\
			0&-1 &0&0&0\\
			0&0& -1&1&0\\
			0&0& 0 &-1&0 \\
			0&0&0&0&1
		\end{smallmatrix} \right]$ via some $P\in \GL(5,\R)$. Setting $f_j=P e_j$ for $1\leq j\leq 5$, it follows from Theorem \ref{theorem: yamada} that $G$ admits a lattice $\Gamma=\pi \Z \ltimes P \Z^5$.
		The form $\tau$ is invariant by $\Gamma$ since  $\exp(-2i(t+\pi k))=\exp(-2it)$. Then $K_{(\Gamma\backslash G,J)}$ is holomorphically trivial.
		
		We remark that $\Gamma\backslash G$ is not homeomorphic to a nilmanifold since the lattice $\Gamma$ is not nilpotent. Indeed, if $\Gamma_{k}$ denotes the $k$-th term of the lower central series of $\Gamma$, then it is easy to compute $\Gamma_{k}=0\Z \oplus \operatorname{Im}(B-\operatorname{I}_5)^k$ which is not trivial because $B-\operatorname{I}_5$ is not nilpotent. 
	\end{example}
	
	

	\section{Applications to hypercomplex geometry} \label{S: hcpx}
	
	In this last section, we explore the triviality of the canonical bundle of complex manifolds obtained from a hypercomplex Lie group $(G,\{J_1,J_2,J_3\})$, or the corresponding quotients by uniform lattices. More concretely, we will show that if there exists a left invariant trivializing section of $K_{(G,J_\alpha)}$ for some $\alpha=1,2,3$, then any associated compact quotient $(\Gamma\backslash G, J_\alpha)$ has trivial canonical bundle for all $\alpha$, also via an invariant section. However, if the trivializing section of $(\Gamma\backslash G,J_\alpha)$ is not invariant, then $K_{(\Gamma\backslash G,J_\beta)}$ is not necessarily trivial for $\beta\neq\alpha$. Using these results we provide a negative answer to a question by Verbitsky. 
	
	We begin by recalling some facts about hypercomplex manifolds. A hypercomplex structure on $M$ is a triple of complex structures $\{J_1,J_2,J_3\}$ on $M$ which obey the laws of the quaternions:
	\[  J_1J_2=-J_2J_1=J_3.  \]
	In particular, $J_\alpha J_\beta =- J_\beta J_\alpha=J_\gamma$ for any cyclic permutation $(\alpha,\beta,\gamma)$ of $(1,2,3)$. 
	
	It follows that $M$ carries a 2-sphere of complex structures. Indeed, given $a=(a_1,a_2,a_3)\in \mathbb{S}^2$,
	\begin{equation}\label{eq:J-sphere}
		J_a:=a_1J_1+a_2J_2+a_3J_3
	\end{equation}
	is a complex structure on $M$. Moreover, for any $p\in M$, the tangent space $T_p M$ has an $\mathbb{H}$-module structure, where $\mathbb{H}$ denotes the quaternions. In particular $\dim_\R M=4n$, $n\in \N$. 
	
	Any hypercomplex structure $\{J_\alpha\}$ on $M$ determines a unique torsion-free connection $\nabla^{\mathcal{O}}$, called the \textit{Obata connection} (see \cite{Ob}), which satisfies $\nabla^{\mathcal{O}} J_\alpha=0$ for all $\alpha$. It was shown in \cite{Sol} that an expression for this connection is given by:
	\begin{equation*}\label{eq:Obata}
		\nabla^\O_X Y=\tfrac12 \left([X,Y]+J_1[J_1 X,Y]-J_2[X,J_2 Y]+J_3 [J_1 X,J_2 Y]\right), \quad X,Y\in\X(M).
	\end{equation*}
	
	Given the hypercomplex structure $\{J_\alpha\}$ on $\R^{4n}$ induced by the quaternions, we denote by
	\[ \GL(n,\H):=\{T\in \GL(4n,\R): TJ_\alpha=J_\alpha T\; \text{for all}\, \alpha\},\] the quaternionic general linear group, with corresponding Lie algebra
	\[ \mathfrak{gl}(n,\H):=\{T\in \gl(4n,\R): TJ_\alpha=J_\alpha T\; \text{for all}\, \alpha\}.\]
	Since $\nabla^\O J_\alpha=0$ for all $\alpha$, the holonomy group of the Obata connection, $\operatorname{Hol}(\nabla^{\mathcal{O}})$, is contained in $\GL(n,\H)$. 
	A hypercomplex manifold $(M^{4n},\{J_\alpha\})$ is said to be an $\SL(n,\H)$-manifold if $\operatorname{Hol}(\nabla^\O)\subset \SL(n,\H)$, where $\SL(n,\H)=[\GL(n,\H),\GL(n,\H)]$ is the commutator subgroup of $\GL(n,\H)$. These manifolds have been actively studied (see for instance \cite{GeTa,GP,GLV,IP,LW,LW1}).
	
	\medskip
	
	We will consider now left invariant hypercomplex structures on Lie groups, which are given equivalently by hypercomplex structures on Lie algebras, as usual. The corresponding Obata connection is also left invariant and it can be determined by its action on left invariant vector fields, that is, on the Lie algebra.
	
	As an application of Theorem \ref{theorem: inv} we show that if $(G^{4n},J_\alpha)$ admits a non-vanishing left invariant closed $(2n,0)$-form for some $\alpha=1,2,3$, then $(G^{4n},J_a)$ (with $J_a$ as in \eqref{eq:J-sphere}) has a non-vanishing left invariant closed $(2n,0)$-form, for all $a\in \mathbb{S}^2$. 
	
	\begin{theorem}\label{theorem: hcpx}
		Let $\{J_1,J_2,J_3\}$ be a hypercomplex structure on the $4n$-dimensional Lie algebra $\g$. If $J_\alpha$ admits a non-vanishing closed $(2n,0)$-form for some $\alpha=1,2,3$, then $J_a$ admits a non-vanishing closed $(2n,0)$-form for any $a\in \mathbb{S}^2$, with $J_a$ given by \eqref{eq:J-sphere}.
	\end{theorem}
	
	\begin{proof}
		Let $(\alpha,\beta,\gamma)$ a cyclic permutation of $(1,2,3)$ with $J_\alpha$ satisfying the conditions in the statement. Then, due to the vanishing of the Nijenhuis tensor $N_{J_\gamma}$, for any $x,y\in\g$ we get 
		\[ J_\gamma [x,y]=[J_\gamma x,y]+[x,J_\gamma y]+J_\gamma [J_\gamma x,J_\gamma y].\] 
		Since $J_\gamma=J_\alpha J_\beta$, applying $-J_\alpha$ in both sides of this equality we have
		\[
		J_\beta[x,y]  =-J_\alpha [J_\gamma x,y]-J_\alpha [x,J_\gamma y] +J_\beta [J_\gamma x,J_\gamma y],
		\]
		which implies 
		\[ \Tr(J_\beta \ad(x))=-\Tr(J_\alpha \ad(J_\gamma x))-\Tr(J_\alpha \ad(x) J_\gamma) + \Tr(J_\beta \ad(J_\gamma x) J_\gamma). \] 
		Using that $\Tr(AB)=\Tr(BA)$ and $\Tr(J_\alpha \ad (x))=\Tr(\ad(J_\alpha x))$ for all $x\in\g$ due to Theorem \ref{theorem: inv} (since the Koszul 1-form $\psi_\alpha$ vanishes) we arrive at
		\begin{align*}
			\Tr(J_\beta \ad(x)) & = -\Tr(\ad(J_\alpha J_\gamma x))-\Tr(J_\gamma J_\alpha \ad(x))+\Tr(J_\gamma J_\beta \ad(J_\gamma x))\\
			& = \Tr(\ad(J_\beta x))-\Tr(J_\beta \ad(x))-\Tr(J_\alpha \ad(J_\gamma x))\\
			& = \Tr(\ad(J_\beta x))-\Tr(J_\beta \ad(x))+\Tr(\ad(J_\beta x)),
		\end{align*} 
		which implies  $\Tr(J_\beta \ad(x)) = \Tr(\ad(J_\beta x)).$ 
		
		The same computation with $\alpha$ replaced by $\beta$ shows that the same condition holds for $J_\gamma$. It follows that 
		$ \Tr(J_a \ad(x)) = \Tr(\ad(J_a x))$
		for any $a\in \mathbb{S}^2$. Therefore, the corresponding Koszul 1-form $\psi_a$ vanishes and according to Theorem \ref{theorem: inv} the proof is complete.
	\end{proof}
	
	\begin{corollary}
		Let $\{J_1,J_2,J_3\}$ be a hypercomplex structure on the simply connected Lie group $G$ and let $\Gamma$ be a uniform lattice on $G$. If there exists an invariant trivializing section of $K_{(\Gamma\backslash G,J_\alpha)}$ for some $\alpha=1,2,3$, then $(\Gamma\backslash G,J_a)$ has trivial canonical bundle for any $a\in\mathbb{S}^2$.
	\end{corollary}
	
	\medskip
	
	In \cite{Ver}, Verbitsky proves that if $(M,\{J_\alpha\})$ is an $\SL(n,\H)$-manifold then the complex manifold $(M,J_\alpha)$ has trivial canonical bundle for all $\alpha$. Then he poses the following question:
	
	\smallskip
	
	\textbf{Question \cite{Ver}:}  Let $(M,\{J_\alpha\})$ be a compact hypercomplex manifold. If the complex manifold $(M,J_1)$ has trivial canonical bundle, does it follow that $M$ is an $\SL(n,\H)$-manifold? 
	
	\smallskip
	
	In certain cases there is an affirmative answer to this question, for instance when $(M,\{J_\alpha\})$ admits a hyperKähler with torsion metric (\cite[Theorem 2.3]{Ver}) or when $M$ is a hypercomplex nilmanifold 
	\cite[Corollary 3.3]{BDV}. In the latter case, the key fact is that every complex nilmanifold has trivial canonical bundle via an invariant trivializing section. Using the same arguments, in \cite{GeTa} it is proved that if a hypercomplex solvmanifold $(\Gamma\backslash G, \{J_\alpha\})$ admits an invariant trivializing section of $K_{(\Gamma\backslash G,J_\alpha)}$ for some $\alpha$ then $\Gamma\backslash G$ is an $\SL(n,\H)$-manifold. 
	
	\smallskip
	
	We exhibit next an hypercomplex solvmanifold $(\Gamma\backslash G,\{J_\alpha\})$ such that $K_{(\Gamma\backslash G,J_1)}$ is trivial but the solvmanifold is not an $\SL(n,\H)$-manifold. 
	
	\begin{example}\label{ex: hcpx}
		Let $\g=\text{span}\{e_1,\ldots,e_4\}$ be the 4-dimensional unimodular 
		Lie algebra given by \[  [e_2,e_3]=e_1, \ \quad [e_2,e_4]=e_2, \quad [e_3,e_4]=-e_3.\]
		It is easily verified that the almost complex structure $J$ defined by $Je_1=e_2$ and $Je_3=e_4$ is integrable. Note that $\g_+=\text{span}\{e_1,e_3\}$ and $\g_-:=J\g_+=\text{span}\{e_2,e_4\}$ are subalgebras of $\g$. Then, according to \cite{AS}, the Lie algebra $\hat{\g}:=(\g_\C)_\R$ admits a hypercomplex structure $\{J_1,J_2,J_3\}$. Indeed, with respect to the decomposition $\hat{\g}=\g\oplus i\g$, these complex structures are given by
		\begin{gather*}
			J_1(x+iy)=\begin{cases} 
				i(x+iy), & x,y\in\g_+, \\
				-i(x+iy),&  x,y\in\g_-,
			\end{cases}\\
			J_2(x+iy)=J x+i J y, \quad x,y\in \g, 
		\end{gather*}
		and $J_3=J_1J_2$. Let us write them down explicitly. Relabelling the basis $\{e_1,\ldots,e_4,i e_1,\ldots,ie_4\}$ as $\{e_1,\ldots,e_8\}$ we have that $\{J_1,J_2,J_3\}$
		are given by
		\begin{align*} 
			J_1 e_1=e_5, &\quad J_1 e_2=-e_6, \quad J_1 e_3=e_7, \quad J_1 e_4=-e_8,\\
			J_2 e_1=e_2, &\quad J_2 e_3=e_4, \quad J_2 e_5=e_6, \quad J_2 e_7=e_8,\\
			J_3 e_1=-e_6, &\quad J_3 e_2=-e_5, \quad J_3 e_3=-e_8, \quad J_3 e_4=-e_7.
		\end{align*}
		With respect to this basis the Lie brackets of $\hat{\g}$ are
		\begin{align*}
			[e_2,e_4]=e_2, \; [e_3,e_4]=-e_3, \; [e_4,e_6]=-e_6,\; [e_4, e_7]=e_7,\\
			[e_2,e_8]=e_6, \; [e_3,e_8]=-e_7, \; [e_6,e_8]=-e_2, \; [e_7,e_8]=e_3, \\
			[e_2,e_3]=e_1, \; [e_2,e_7]=e_5, \; [e_3,e_6]=-e_5, \; [e_6,e_7]=-e_1.
		\end{align*}
		If we denote $\psi_{\alpha}(x):=\Tr(J_\alpha \ad x)$ for $\alpha=1,2,3$, then 
		\[ \psi_1=-4 e^8, \quad \psi_2=-4 e^3, \quad \psi_3=-4 e^7,\]
		where $\{e^j\}_{j=1}^8$ is the dual basis of $\{e_j\}_{j=1}^8$. Since $\psi_{\alpha}\neq 0$, we have that $(\hat{\g},J_\alpha)$ does not admit a nonzero closed $(4,0)$-form, for any $\alpha$.
		
		Moreover, note that $\psi_1([\hat{\g},\hat{\g}])=0$ but $\psi_2([\hat{\g},\hat{\g}])\neq 0$ and $\psi_3([\hat{\g},\hat{\g}])\neq 0$. According to Theorem \ref{theorem: obstruction}, for any lattice $\Lambda\subset \hat{G}$, where $\hat{G}$ is the simply connected Lie group associated to $\hat{\g}$, the compact complex manifold $(\Lambda\backslash \hat{G},J_{\alpha})$ has non-trivial canonical bundle for $\alpha=2,3$. Nevertheless we show next that there exist lattices $\Gamma_m\subset \hat{G}$  such that the corresponding complex solvmanifolds $(\Gamma_m\backslash \hat{G},J_1)$ do have trivial canonical bundle.
		
		We show first that $\hat{G}$ admits a lattice $\Gamma_m$ for any $m\in\N, \, m\geq 3$. Indeed, we may write $\hat{\g}=(\R e_8 \oplus \R e_4)\ltimes \n$, where the nilradical $\n$ is spanned by $\{e_1,e_2,e_3,e_5,e_6,e_7\}$. We compute
		\begin{align*} A:=\exp(\pi \ad e_8|_{\n})&=\operatorname{diag}(1,-1,-1,1,-1,-1), \\
			B_m:=\exp(t_m \ad e_4|_{\n})&=\operatorname{diag}(1,\alpha_m^{-1},\alpha_m,1,\alpha_m^{-1},\alpha_m),
		\end{align*}
		where $\alpha_m=\frac{m+\sqrt{m^2-4}}{2}$, and $t_m=\log \alpha_m$.
		Setting 
		\[P_m=\left[\begin{array}{ccc} 
			1&0&0\\ 
			0&1&\alpha_m^{-1}  \vspace{0.1cm}\\ 
			0 & \frac{1}{\alpha_m-\alpha_m^{-1}}& \frac{\alpha_m}{\alpha_m-\alpha_m^{-1}}
		\end{array}\right]^{\oplus 2},\] we obtain $P_m^{-1} B_m P_m=\left[\begin{smallmatrix}1&0&0\\0&0&-1\\0&1&m \end{smallmatrix}\right]^{\oplus 2}$ and $P_m^{-1} A P_m=A$ for any $m$. If we define $f_j=P_m e_j$ for $j=1,2,3,5,6,7$ then it is easy to verify that $[f_k,f_\l]=[e_k,e_\l]$ for $k,\l\in\{1,2,3,5,6,7\}$. Therefore, $\{f_1,f_2,f_3,f_5,f_6,f_7\}$ is a rational basis of $\n$ in which $A$ and $B_m$ are expressed as unimodular integer matrices. According to Theorem \ref{theorem: yamada}, the semidirect product \[\Gamma_m=(\pi \Z \oplus t_m \Z)\ltimes \exp^N (\text{span}_\Z\{f_1,f_2,f_3,f_5,f_6,f_7\})\] is a lattice in $\hat{G}=\R^2\ltimes N$, where $N$ is the nilradical of $\hat{G}$.
		
		Let $\sigma_1:=(e^1+ie^5) \wedge (e^2-ie^6) \wedge (e^3+ie^7) \wedge (e^4-ie^8)$ which is a nonzero $(4,0)$-form with respect to $J_1$. It follows from Proposition \ref{proposition: tau explicit} and $-\frac{1}{2} \Tr(J_1\ad e_8)=2$ that $\tau_1:=\exp(2i x_8) \sigma_1$ is a nonzero closed $(4,0)$-form on $\hat{G}$ with respect to $J_1$, where $\mathbf{x}:=(x_8,x_4,x_1,x_2,x_3,x_5,x_6,x_7)$ are the real coordinates of $\hat{G}$. It follows from $\exp(2 i (x_8+\pi k))=\exp(2 i x_8)$ that $f(\mathbf{x})=\exp(2 i x_8)$ is invariant by the action of $\Gamma_m$ so there is an induced smooth function $\hat{f}:\Gamma_m\backslash \hat{G}\to\C$ such that the $(4,0)$-form $\hat{\tau}_1=\hat{f} \hat{\sigma}_1$ is a trivializing section of $(\Gamma_m\backslash \hat{G},J_1)$. In particular, $(\Gamma_m\backslash \hat{G}, J_1)$ has trivial canonical bundle.
		
		If $\operatorname{Hol}(\nabla^\O)$ were contained in $\SL(n,\H)$, then the canonical bundle of $(\Gamma_m\backslash \hat{G},J_\alpha)$ would be trivial for all $\alpha$ but we have shown that this is not the case for $\alpha=2$ and $\alpha=3$. Therefore, this example provides a negative answer to Verbitsky's question.
	\end{example}
	
	\begin{remark}
		Verbitsky's question remains open for a hypercomplex manifold $(M,\{J_\alpha\})$ such that $(M,J_{\alpha})$ has trivial canonical bundle for all $\alpha$.
	\end{remark}
	
	\medskip

\textbf{Conflict of interest.} 	The authors have no conflict of interest to declare that are relevant to this article.

	\ 
	
\end{document}